\documentclass[12pt]{article}
\usepackage[margin=1in]{geometry}

\usepackage{graphicx}
\usepackage{float}

\usepackage{amssymb}
\usepackage{amsmath}
\usepackage{amsthm}

\newtheorem{theorem}{Theorem}
\newtheorem{lemma}[theorem]{Lemma}

\newtheorem{proposition}[theorem]{Proposition}
\theoremstyle{definition}
\newtheorem{definition}[theorem]{Definition}
\newtheorem{remark}[theorem]{Remark}
\newtheorem*{uremark}{Remark}

\usepackage{bm}
\usepackage[authoryear]{natbib}
\usepackage[hidelinks]{hyperref}
\usepackage{comment}
\usepackage{color}
\usepackage{tikz}

\newcommand{\Z}{\mathbb{Z}}

\newcommand{\N}{\mathbb{N}}
\newcommand{\R}{\mathbb{R}}
\newcommand{\T}{\mathbb{T}}

\renewcommand{\L}{\mathcal{L}}

\renewcommand{\a}{\alpha}

\newcommand{\bma}{\bm{\alpha}}
\renewcommand{\d}{\delta}

\newcommand{\e}{\varepsilon}
\newcommand{\signvec}{\boldsymbol{\epsilon}}
\renewcommand{\i}{\mathbf{i}}

\newcommand{\mbf}{\mathbf}

\newcommand{\SlProjtwo}{S(\mbf{l}, \signvec;N)}
\newcommand{\ElProjtwo}{E (\mbf{l},\signvec;N)}

\newcommand{\productinDProjtwobar}{\prod_{j=1}^{d}\l(\frac{1-e^{-2\pi \i N(n_1\alpha_j -n_{j+1})}}{2\pi (n_1\alpha_j -n_{j+1})}\r)\(\fc{\sin 2\pi (n_1\a_j -n_{j+1})}{2\pi (n_1\a_j-n_{j+1})}\)^2}
\newcommand{\convergentconditionforphi}{Let $\varphi(n)$ be a positive nondecreasing function of $n$ with $\sum_{n=1}^{\infty}\fc{1}{\varphi(n)}<\infty$}

\newcommand{\beq}{\begin{equation}}
\newcommand{\eeq}{\end{equation}}

\newcommand{\fc}{\frac}
\newcommand{\nid}{\noindent}

\renewcommand{\l}{\left}
\renewcommand{\r}{\right}
\newcommand{\sums}{\sum\limits}

\newcommand{\ba}{\begin{array}}
\newcommand{\ea}{\end{array}}
\newcommand{\bal}{\begin{aligned}}
\newcommand{\eal}{\end{aligned}}

\renewcommand{\(}{\left(}
\renewcommand{\)}{\right)}

\title{Almost Sure Bounds for Discrepancies of Linear Forms on the Circle}
\author{Hao Wu\thanks{Institute of Analysis and Number Theory, TU Graz, Austria. Email: \texttt{haowu.nankai@gmail.com}.}}
\date{}

\begin{document}

\maketitle

\begin{abstract}
As a generalization of irrational rotations and a dual case of higher-dimensional Kronecker sequences, we study the discrepancy of sequences of linear forms on the circle. Given irrationals $\alpha_1,\dots,\alpha_d$, consider the set of $N^d$ points $\{k_1\alpha_1+\cdots+k_d\alpha_d \mod 1 : 1\le k_j\le N\}$. We prove that for a full-measure set of vectors $(\alpha_1,\dots,\alpha_d)\in\mathbb{R}^d$, the maximal discrepancy of these points relative to intervals in $[0,1)$ has the optimal principal order $(\log N)^d$, up to powers of $\log\log N$. This result provides a nearly sharp dual analogue, in the setting of linear forms, to Beck's celebrated theorem on multidimensional Kronecker sequences (Ann. of Math., 1994). The proof combines Fourier analysis, metric multiplicative Diophantine estimates, and a duality argument which reduces certain lattice-counting errors to Beck's discrepancy theorem.
\end{abstract}


\section{Introduction}

Let $\alpha\in \R$ be irrational, and denote $\T=\R/\Z\cong [0,1)$ the $1$-dimensional torus. By Weyl's criterion (see \cite{Drmota}, Section 1.2.1), the sequence $\{n\a\mod 1\}_{1\le n\le N}$ is equidistributed on $\T$ as $N\to \infty$; that is, for any interval $[a,b)\subset \T$,
$$
\frac{1}{N}\sum_{n=1}^{N}\chi_{[a,b)}(n\a)\rightarrow b-a, \quad N\rightarrow \infty.
$$

To measure the rate of convergence, we introduce the discrepancy function, defined as the difference between the actual number of visits to $[a,b)$ up to time $N$ and the expected number of visits $N\cdot (b-a)$:
$$
D_{[a,b)}(\a;N)=\sum_{n=1}^{N}\chi_{[a,b)}(n\a\mod 1)-N(b-a).
$$

In the 1920s, \cite{Khintchine} proved that for any $\e>0$, for almost every irrational ${\alpha\in \R}$, the maximal discrepancy relative to all intervals in $[0,1)$:
$$
\Delta(\a;N)=\max_{x\in(0,1]} |D_{[0,x)}(\a;N)|
$$
satisfies
$$
\Delta(\a;N)\ll_{\a,\e}\log N(\log \log N)^{1+\e},
$$
and also exceeds $\log N\cdot \log\log N$ infinitely often. Khintchine's proof relies heavily on the classical metrical theory of continued fractions.

Due to the lack of an effective continued fraction algorithm in higher dimensions, the higher-dimensional analogue of Khintchine's theorem remained open for a long time. \cite{schmidt1964metrical} proved that $\Delta(\bma;N)\ll (\log N)^{d+1+\e}$ for almost every $\bma\in\R^d$, using the Erdős--Turán--Koksma inequality.

In a celebrated paper, \cite{Beck} removed the extra $\log N$ factor by combining Fourier analysis, the second moment method, and combinatorics. More precisely, Beck proved the following multidimensional analogue of Khintchine's theorem.

Let $d\ge 2$, let $\bm{\a}=(\a_1,\dots,\a_d)\in \R^d$ be the translation vector, and let
$$
B(\mbf{x})=[0,x_1)\times \dots \times [0,x_d)\subset [0,1)^d.
$$
Define
$$
D(\bm{\a},\mbf{x};m)=\sum_{1\le n\le m} \chi_{B(\mbf{x})}(n\bm{\a})-m\text{Vol}(B(\mbf{x}))
$$
and
\beq\label{Beck maximal discrepancy}
\Delta(\bm{\a};N)= \max_{\mbox{\scriptsize$\begin{array}{c}
1\le m \le N \\
\mbf{x}\in [0,1]^d\\\end{array}$}} |D(\bm{\a},\mbf{x};m)|.
\eeq
Then, for every positive nondecreasing function $\varphi(n)$, for almost every $\bm{\a}\in \R^d$,
\beq\label{Thm of Beck}
\Delta(\bm{\a};N)\ll_{\bma,\varphi} (\log N)^d\cdot \varphi(\log \log N)
\Longleftrightarrow
\sum_{n=1}^\infty \frac{1}{\varphi(n)}<\infty.
\eeq



In this paper, we study the discrepancy of the sequence of linear forms
\[
\left\{
\sum_{1\le j\le d} k_j\a_j \mod 1
\right\}_{\substack{1\le k_j\le N \\ 1\le j\le d}}
\]
relative to intervals in $[0,1)$.

The discrepancy problem for linear forms may be viewed as a dual counterpart of the classical Kronecker discrepancy problem. 
After Fourier expansion, the principal contribution is governed by the multiplicative resonant divisors $n\prod_{j=1}^{d}\|n\alpha_j\|$, which are products of the coordinates of the associated Kronecker sequence.
A priori, such multiplicative resonances could lead to an upper bound of order $(\log N)^{d+1}$, analogous to the additive resonances $\prod_j^d n_j \cdot\|n_1\alpha_1+\cdots+n_d\alpha_d\|$ in the pre-Beck situation.


To state our main theorem, define the discrepancy function of linear forms by
\beq\label{original discrepancy function}
D(\bma,x;N)
=
\sum \limits_{\substack{1\le k_j\le N \\ 1\le j\le d}}
\chi_{[0,x)}
\left(
\sum\limits_{1\le j\le d} k_j\a_j \mod 1
\right)
-
N^d x,
\eeq
and the maximal discrepancy by
$$
\Delta(\bma; N)=\max_{0<x\le 1} |D(\bma,x;N)|.
$$

Informally, our main theorem shows that for almost every $\bm{\alpha}\in\R^d$, the maximal discrepancy of linear forms has essentially optimal order $(\log N)^d$, up to powers of $\log\log N$. The precise statement is the following.
\begin{theorem}\label{Proj two main result}
Let $\varphi(n)$ be a positive nondecreasing function of $n$. Then, for $d\ge 2$ and for almost every $\bm{\a}\in \R^d$, we have
\beq\label{Proj two convergent part of the theorem}
\sum_{n=1}^\infty \fc{1}{\varphi(n)}<\infty
\Rightarrow
\Delta(\bm{\a}; N)\ll_{\bm{\a},\varphi} (\log N)^d\cdot \varphi^3(\log \log N),
\eeq
\beq\label{Proj two divergent part of the theorem}
\sum_{n=1}^\infty \fc{1}{\varphi(n)}=\infty
\Rightarrow
\Delta(\bm{\a}; N)\gg_{\bm{\a},\varphi} (\log N)^d\cdot \varphi(\log \log N)
\quad\text{infinitely often in }N.
\eeq
\end{theorem}
\begin{uremark}

The main novelty of Theorem \ref{Proj two main result} is that we preserve the optimal principal order $(\log N)^d$ in the presence of the aforementioned multiplicative resonant divisors $n\prod_{j=1}^{d}\|n\a_j\|$.

The hope of achieving a full analogue of the $0$--$1$ law in \eqref{Thm of Beck} relies on a sufficiently uniform metrical theorem in multiplicative Diophantine approximation, as also sought by \cite{bjorklund2023uniformmetricaltheoremmultiplicative} in dimension $d=2$, but for significantly smaller lattice counting domains and in arbitrary dimensions.

Instead, we use a first moment estimate together with a double exponential dyadic decomposition (Proposition \ref{Proj two control for sum of small divisors}), and exploit the duality between linear forms and Kronecker sequences (Lemma \ref{Proj two number of elements in S for N^d/4} and \ref{Proj two d_N error estimation for the number of elements}). This dual use of Beck's theorem supplies the required uniform counting estimate for the dyadic cells appearing in the linear-form problem. The argument therefore does not rely on the second moment method used in \cite{Beck}.

\end{uremark}

\textbf{Motivation.}

The sequence of linear forms arises naturally at the intersection of number theory, ergodic theory, and mathematical physics.

From the dynamical and ergodic-theoretic point of view, the shifted family
\[
\left\{
x+\sum_{j=1}^{d} k_j\alpha_j \mod 1
:
1\le k_j\le N
\right\}
\]
is the orbit of the $\mathbb Z^d$-action generated by the circle rotations $x\mapsto x+\alpha_j \mod 1$, observed through interval test functions and along cubic F\o lner sets. Thus Theorem \ref{Proj two main result} may be read as an almost-sure deviation estimate for this higher-rank action. In analogy with the one-dimensional limit laws of \cite{Kesten1,Kesten2} and the recent refinement and local-discrepancy application of \cite{Borda2025}, and with the higher-dimensional work of \cite{dolgopyat2014deviations,dolgopyat2012deviations}, it is natural to ask for limiting distributions after randomizing the starting point $x$ and the parameters $\alpha_j$. The present paper does not address this limit-law question, but the scale in Theorem \ref{Proj two main result} suggests that the natural normalization should be of order $(\log N)^d$.


In mathematical physics, the Berry--Tabor conjecture \cite{BerryTabor1977} predicts Poisson statistics for the energy gaps in generic integrable systems; the multidimensional harmonic oscillator is a notable exception, whose energy level spacings reduce directly to the fractional parts of a linear form. For algebraic frequencies, \cite{HaynesRoeder2022} proved that these spacings depend quasiperiodically on $\log N$. In contrast to the Poisson regime, 
we focus on the macroscopic behavior of the sequence of linear forms for a full-measure set of frequencies. At this macroscopic level, the logarithmic discrepancy scale in \eqref{Proj two convergent part of the theorem} is different from the polynomial scale suggested by independent Poisson counting.


A closely related direction is the higher-dimensional three-gap problem. Using homogeneous dynamics, \cite{HaynesMarklof2020} proved that for almost every $\boldsymbol{\alpha}$ the number of distinct gap lengths in the sequence of linear forms is unbounded, with a close connection to the Littlewood conjecture.

Finally, for Kronecker sequences the connection between discrepancy and multiplicative Diophantine approximation was formalized by \cite{Beck} and was recently made deterministic by \cite{ChowTechnau2024}. In \cite{ChowTechnau2024}, the authors asked whether an analogous deterministic result could also hold for the discrepancy of dual Kronecker sequences. Our result concerns the full-measure setting and may be viewed as a first step in that direction.

\vspace{+0.2cm}

\textbf{Plan of the paper.}

Overall, our approach is inspired by \cite{Beck}; we adapt his methods to treat the multiplicative resonant divisors $n\prod_{j=1}^{d}\|n\a_j\|$. The paper is organized as follows.

In Section \ref{section 1 of poisson summation formula}, we transform the discrepancy function into its Fourier series using Poisson's summation formula. The convergence of the series follows from a roof-like smoothing procedure combined with an induction on the dimension.

In Section \ref{section 2 estimation of the tails}, we estimate the contribution of the tails and highly resonant Fourier modes. We show that the high-frequency modes and large divisors contribute at most $(\log\log N)^d\log\log N$. The highly resonant terms, whose divisors can be as small as $\frac{1}{(\log N)^d\log\log N}$, are shown to contribute at most $(\log N)^d\varphi^3(\log\log N)$.

Sections \ref{Section 3 of cancellation of the main terms.} and \ref{section 5 estimation of small exponentials} treat the "main" part of the discrepancy. In Section \ref{Section 3 of cancellation of the main terms.}, we pair Fourier modes with opposite signs to achieve cancellation, and thus control the contribution from the part with constant numerators. In Section \ref{section 5 estimation of small exponentials} we adapt Beck's \textit{line method} to estimate the remaining oscillatory part, by showing that modes without sufficient cancellation in the exponential part are sparse. 

Finally, in Section \ref{divergent part of the main theorem}, we prove the divergent statement \eqref{Proj two divergent part of the theorem}. Using a local theorem in simultaneous Diophantine approximation in \cite{gallagher1962metric}, we construct a resonant Fourier mode whose contribution is of order $(\log N)^d(\log\log N)$, and then show that this single mode is dominant.

\section{Poisson's summation formula}\label{section 1 of poisson summation formula}

In this section, following \cite{Beck}, we apply Poisson's summation formula to obtain a Fourier series representation of the discrepancy function. The main results are the explicit formula \eqref{Final expression of the smoothed discrepancy} for a smoothed discrepancy function $\bar{D}$ and Proposition \ref{Proj two Prop Dbar}. 

Throughout the remainder of the paper, we work with the series \eqref{Final expression of the smoothed discrepancy} instead of the original discrepancy function \eqref{original discrepancy function}. The formula \eqref{Final expression of the smoothed discrepancy} defines a smoothed version of $D(\boldsymbol{\alpha},x;N)$ (see \eqref{Smoothed discrepancy function}) with improved convergence properties. Proposition \ref{Proj two Prop Dbar} shows that the difference between $D(\boldsymbol{\alpha},x;N)$ and $\bar{D}$ is sufficiently small to justify this replacement.

The convergent upper bound in Theorem \ref{Proj two main result} is proved by induction on the dimension. The base input is the one-dimensional Khintchine theorem, equivalently the case $d=1$ of \eqref{Thm of Beck}. Assuming the convergent statement already known for all dimensions $<d$, Proposition \ref{Proj two Prop Dbar} gives the smoothing comparison in dimension $d$, while Sections \ref{section 2 estimation of the tails}--\ref{section 5 estimation of small exponentials} prove the required bound for the smoothed discrepancy in the same dimension, which closes the induction. The divergent statement in Section \ref{divergent part of the main theorem} is independent of this induction.

\subsection{Notations}
 
For a real number $x$ we write $[x]$ for its integer part, $\{x\} = x - [x]$ for
its fractional part, $\lceil x\rceil=\min\{y\in \Z,\  y\ge x\}$ for the closest integer from above, and 
\[
  \|x\| = \min\bigl\{|x - [x]|,\, |x - ([x]+1)\}
\]
for its distance to the nearest integer. Also, $f(N) \asymp_d g(N)$ means that $f(N) \ll_d g(N)$ and $g(N)\ll_d f(N)$, where the constants depend on $d$, but not on $N$.
\medskip

\subsection{A formal Fourier expansion}

We begin with a heuristic derivation of a \emph{formal} Fourier expansion of the discrepancy function. Ignoring convergence issues for the moment, Poisson summation suggests that
\begin{equation}\label{formal Fourier series}
D(\boldsymbol{\alpha},x;N)
=\mathrm{i}^{\,d+1} \sum_{\mathbf{n}\in \mathbb{Z}^{d+1}\setminus\{0\}} 
\frac{1-e^{2\pi \mathrm{i} n_1 x}}{2\pi n_1}
\prod_{j=1}^{d}
\frac{1-e^{-2\pi \mathrm{i} N(n_1\alpha_j - n_{j+1})}}{2\pi (n_1\alpha_j - n_{j+1})}.
\end{equation}
Here by continuity, the expression $(1-e^{2\pi \mathrm{i} n_1 x})/(2\pi n_1)$ is interpreted as $-\mathrm{i}x$ when $n_1=0$, and similarly each factor 
\[
\frac{1-e^{-2\pi \mathrm{i} N(n_1\alpha_j - n_{j+1})}}{2\pi (n_1\alpha_j - n_{j+1})}
\]
is interpreted as $\mathrm{i}N$ when $n_1\alpha_j - n_{j+1}=0$.

\medskip

To motivate \eqref{formal Fourier series}, observe that the condition
\[
0 \le \left\{\sum_{j=1}^d k_j \alpha_j \right\} < x
\]
is equivalent to the existence of $m \in \mathbb{Z}$ such that
\[
0 \le \sum_{j=1}^d k_j \alpha_j - m < x,
\]
with $1 \le k_j \le N$. This embeds the problem into the lattice
\[
\left\{\left(\sum_{j=1}^d k_j \alpha_j - m, k_1, \dots, k_d\right) : (k_1,\dots,k_d,m)\in \mathbb{Z}^{d+1}\right\} \subset \mathbb{R}^{d+1}.
\]
The condition above corresponds to lattice points lying in the box
\[
B = [0,x) \times (0,N]^d.
\]

Thus the discrepancy can be written as
\[
\sum_{(k_1,\dots,k_d,m)\in \mathbb{Z}^{d+1}} 
\chi_B\!\left(\sum_{j=1}^d k_j \alpha_j - m, k_1,\dots,k_d\right).
\]

Introducing the matrix
\[
\mathbf{A} =
\begin{pmatrix}
\alpha_1 & \cdots & \alpha_d & -1 \\
1 & \cdots & 0 & 0 \\
\vdots & \ddots & \vdots & \vdots \\
0 & \cdots & 1 & 0
\end{pmatrix},
\]
this expression becomes $\chi_B(\mathbf{A}\mathbf{y})$ with $\mathbf{y}\in\mathbb{Z}^{d+1}$. Formally applying Poisson summation to $f(\mathbf{y}) = \chi_B(\mathbf{A}\mathbf{y})$ yields \eqref{formal Fourier series} after evaluating the resulting integrals.

\medskip

The difficulty with this argument is that Poisson summation requires sufficient smoothness and decay of both $f$ and its Fourier transform. Since indicator functions do not satisfy these conditions, we instead follow \cite{Beck} and introduce a smoothed version of the discrepancy.

\subsection{Smoothing via local averaging}

We average $D(\boldsymbol{\alpha},x;N)$ over small perturbations of both the interval $[0,x)$ and the summation range. For $\mathbf{u}=(u_2,\dots,u_{d+1})$, define the oscillated discrepancy
\begin{equation}\label{Oscillated Discrepancy}
D(\boldsymbol{\alpha}; a,b; \mathbf{u},\mathbf{u}+N)
=\sum_{\substack{u_{j+1}<k_j\le N+u_{j+1}}} 
\chi_{[a,b)}\!\left(\sum_{j=1}^d k_j \alpha_j \bmod 1\right)
- N^d(b-a).
\end{equation}

We then define the smoothed discrepancy by
\begin{equation}\label{Smoothed discrepancy function}
\begin{aligned}
\bar{D}(\boldsymbol{\alpha},x;N)
&=\frac{N^d}{2}\left(\frac{1}{2}\right)^d 
\int_{-2/N^d}^{2/N^d} \int_{[-2,2]^d}
\left(1-\frac{N^d}{2}|u_1|\right)
\prod_{j=1}^d \left(1-\frac{|u_{j+1}|}{2}\right) \\
&\quad \times D(\boldsymbol{\alpha}; u_1, x+u_1; \mathbf{u}, \mathbf{u}+N)\, d\mathbf{u}du_1.
\end{aligned}
\end{equation}

Using standard Fejér kernel identities,
\[
\frac{N^d}{2}\int_{-2/N^d}^{2/N^d}
\left(1-\frac{N^d}{2}|u_1|\right)e^{2\pi \mathrm{i} n_1 u_1}\,du_1
=\left(\frac{\sin\left(2\pi n_1/N^d\right)}{2\pi n_1/N^d}\right)^2,
\]
and for each $1\le j\le d$,
\[
\frac{1}{2}\int_{-2}^{2}
\left(1-\frac{|u_{j+1}|}{2}\right)e^{-2\pi \mathrm{i} (n_1\alpha_j-n_{j+1})u_{j+1}}\,du_{j+1}
=\left(\frac{\sin 2\pi (n_1\alpha_j-n_{j+1})}{2\pi (n_1\alpha_j-n_{j+1})}\right)^2.
\]
The zero Fourier mode is exactly cancelled by the volume term in the definition of the discrepancy. Hence the smoothed discrepancy has the Fourier expansion
\begin{equation}\label{Final expression of the smoothed discrepancy}
\begin{aligned}
\bar{D}(\boldsymbol{\alpha},x;N)
=\mathrm{i}^{\,d+1} \sum_{\mathbf{n}\in \mathbb{Z}^{d+1}\setminus\{0\}}
&\frac{1-e^{2\pi \mathrm{i} n_1 x}}{2\pi n_1}
\left(\frac{\sin\left(2\pi n_1/N^d\right)}{2\pi n_1/N^d}\right)^2 \\
&\cdot \prod_{j=1}^{d}\left(\frac{1-e^{-2\pi \mathrm{i} N(n_1\alpha_j -n_{j+1})}}{2\pi (n_1\alpha_j -n_{j+1})}\right)
\left(\frac{\sin 2\pi (n_1\alpha_j -n_{j+1})}{2\pi (n_1\alpha_j-n_{j+1})}\right)^2.
\end{aligned}
\end{equation}
The smoothing ensures that the Fourier coefficients decay rapidly (essentially like $\prod n_i^{-3}$), so Poisson summation can now be applied. In particular, we will show that the series \eqref{Final expression of the smoothed discrepancy} converges absolutely in Proposition \ref{Proj two D_1 - D convergence of the poisson summation}, which justifies the Poisson summation fully.
\medskip

\subsection{Comparison with the original discrepancy}

Finally, we justify replacing the original discrepancy by its smoothed version. We show that the smoothing procedure introduces only a small error.

\begin{proposition}\label{Proj two Prop Dbar}
For $d\ge 2$, and every $\varepsilon>0$, there exists a full-measure set of $\boldsymbol{\alpha}\in[0,1)^d$ such that
\begin{equation}\label{Difference between the smooth and original discrepancy}
|\bar{D}(\boldsymbol{\alpha},x;N)-D(\boldsymbol{\alpha},x;N)|
\ll_{\boldsymbol{\alpha},\varepsilon} (\log N)^{d-1+\varepsilon}.
\end{equation}
\end{proposition}

\begin{proof}
It suffices to prove that, for every $u_1$ and $\boldsymbol u$ in the integration range of \eqref{Smoothed discrepancy function},
\begin{equation}\label{prop2-target-comparison}
|D(\boldsymbol{\alpha};u_1,x+u_1;\boldsymbol u,\boldsymbol u+N)-D(\boldsymbol{\alpha},x;N)|
\ll_{\boldsymbol{\alpha},\varepsilon}(\log N)^{d-1+\varepsilon},
\end{equation}
with a constant independent of $x,u_1$ and $\boldsymbol u$. The estimate for $\bar D-D$ then follows by averaging. We use the circle convention for intervals; if an interval crosses $0$, it is split into at most two ordinary intervals, which changes only the number of boundary pieces by an absolute factor.

We prove the comparison estimate by induction on the dimension. For dimension $d=2$ we use Khintchine's or Beck's theorem \eqref{Thm of Beck} in dimension $1$. For dimension $d$, the induction hypothesis is the convergent upper bound in Theorem \ref{Proj two main result} in all dimensions strictly smaller than $d$. Since an arbitrary interval on the circle can be written, up to finite endpoint errors of size $O_d(1)$, as a difference of two anchored intervals, this hypothesis applies uniformly to translated intervals as well.

The difference arises from two sources:

\begin{itemize}
    \item[(i)] Points of the original sequence whose fractional parts lie in the boundary intervals 
    \[
    [-2/N^d,0) \quad \text{and} \quad [x, x+2/N^d),
    \]
    whose total contribution we denote by $D_1$:
    \[
      D_1=D(\boldsymbol{\alpha}; -2/N^d, 0; \mathbf{u}, \mathbf{u}+N)+D(\boldsymbol{\alpha}; 0, 2/N^d; \mathbf{u}, \mathbf{u}+N)
      ;\]
    
    \item[(ii)] Points coming from indices outside the range $\{1,\dots,N\}$ due to the perturbation $\mathbf{u}$. We denote by $D_2$ the remaining error. 
\end{itemize}

These two contributions may overlap, but this does not affect the estimates. Figure~\ref{fig:error-terms} is an illustration of the two types of error.

\begin{figure}[!htbp]
\centering
\resizebox{0.98\textwidth}{!}{%

\begin{tikzpicture}[
  font=\small,
  axis/.style={->, >=stealth, thick},
  fuzzy/.style={red!18},
  mainbox/.style={blue!12},
  boundary/.style={red!18},
  note/.style={align=left, font=\scriptsize, text width=4.6cm},
  label/.style={font=\scriptsize},
  title/.style={font=\bfseries\small}
]

\begin{scope}[xshift=0cm,yshift=0cm]
  \node[title, anchor=west] at (0,5.10) {(a) Origin of $D_1$: value space};
  \node[note, text width=7.25cm, anchor=west] at (0,4.35)
  {$D_1$ counts values of $\sum_j k_j\alpha_j\pmod 1$ that fall into the two short boundary intervals created by smoothing the target interval.};

  \draw[axis] (0.10,2.45) -- (7.25,2.45) node[below=3pt,label] {$1$};

  \draw[thick] (1.15,2.25) -- (1.15,2.65);

  \draw[thick] (4.75,2.25) -- (4.75,2.65);

  \fill[blue!10] (1.15,2.20) rectangle (4.75,2.70);
  \draw[blue!70!black, line width=1.5pt] (1.15,2.45) -- (4.75,2.45);
  \node[blue!70!black, label, above=12pt] at (2.95,2.55)
    {original target interval $[0,x)$};

  \fill[fuzzy] (0.35,2.20) rectangle (1.15,2.70);
  \fill[fuzzy] (4.75,2.20) rectangle (5.55,2.70);
  \draw[red!80!black, dashed] (0.35,2.20) rectangle (1.15,2.70);
  \draw[red!80!black, dashed] (4.75,2.20) rectangle (5.55,2.70);

  \filldraw[red!85!black] (0.75,2.45) circle (1.8pt);
  \filldraw[red!85!black] (5.15,2.45) circle (1.8pt);

  \draw[red!80!black, dashed] (0.75,2.20) -- (0.75,2.02);
  \node[red!80!black, label, align=center, text width=2.55cm] at (0.75,1.78)
    {$[-2/N^d,0)$\\$D_1$ point};

  \draw[red!80!black, dashed] (5.15,2.20) -- (5.15,2.02);
  \node[red!80!black, label, align=center, text width=2.70cm] at (5.15,1.78)
    {$[x,x+2/N^d)$\\$D_1$ point};

  \foreach \x in {1.85,3.05,6.20} {
    \filldraw[black] (\x,2.45) circle (1.35pt);
  }
\end{scope}

\begin{scope}[xshift=8.25cm,yshift=0cm]
  \node[title, anchor=west] at (0,5.10) {(b) Origin of $D_2$: index space, $d=2$};
  \node[note, text width=7.2cm, anchor=west] at (0,4.55)
  {$D_2$ comes from indices outside the main summation box but inside the perturbed box.};

  \begin{scope}[xshift=0.15cm,yshift=0.20cm,scale=0.48]
    \draw[axis] (-0.6,-0.6) -- (7.9,-0.6) node[right,label] {$k_1$};
    \draw[axis] (-0.6,-0.6) -- (-0.6,7.9) node[above,label] {$k_2$};

    \fill[boundary] (0,0) rectangle (7,7);
    \draw[red!80!black, dashed, thick] (0,0) rectangle (7,7);
    \fill[mainbox] (1.35,1.35) rectangle (5.65,5.65);
    \draw[blue!70!black, thick] (1.35,1.35) rectangle (5.65,5.65);

    \draw[dashed] (1.35,0) -- (1.35,7);
    \draw[dashed] (5.65,0) -- (5.65,7);
    \draw[dashed] (0,1.35) -- (7,1.35);
    \draw[dashed] (0,5.65) -- (7,5.65);

    \node[blue!70!black, align=center, font=\scriptsize\bfseries] at (3.5,3.5)
      {main range\\$1\le k_1,k_2\le N$};

    \foreach \x/\lab/\col in {0/{-2}/red!80!black,1.35/{1}/black,5.65/{N}/black,7/{N+2}/red!80!black}
      {\draw (\x,-0.76)--(\x,-0.46) node[below=2pt, label, text=\col] {$\lab$};}
    \foreach \y/\lab/\col in {0/{-2}/red!80!black,1.35/{1}/black,5.65/{N}/black,7/{N+2}/red!80!black}
      {\draw (-0.76,\y)--(-0.46,\y) node[left=2pt, label, text=\col] {$\lab$};}

    \draw[red!80!black, ->, >=stealth, thick] (7.15,4.55) -- (8.55,4.55);
    \draw[red!80!black, ->, >=stealth, thick] (4.55,0.45) -- (8.55,1.25);
  \end{scope}

  \node[red!80!black, align=center, font=\scriptsize\bfseries] at (5.75,2.95)
    {$D_2$ boundary region};
  \node[note, text width=4.25cm, anchor=west] at (4.55,1.85)
    {At least one $k_j$ is fixed near a boundary face, for example\\$k_1\in\{-1,0,N+1,N+2\}$.\\Hence the estimate reduces to dimension $d-1$.};
\end{scope}

\end{tikzpicture}
}
\caption{Origins of the two error terms.}
\label{fig:error-terms}
\end{figure}

First consider $D_1$ due to the change of the target interval. By a standard Borel--Cantelli argument below, for almost every $\boldsymbol{\alpha}$ there is a constant $c(\boldsymbol{\alpha},\varepsilon)>0$ such that, for every nonzero $\boldsymbol h\in\mathbb Z^d$ with $|h_j|\le N+4$,
\begin{equation}\label{prop2-gap-estimate}
 \bigl\|\sum_{j=1}^d h_j\alpha_j\bigr\|
 \ge \frac{c(\boldsymbol{\alpha},\varepsilon)}{N^d(\log N)^{1+\varepsilon}}.
\end{equation}
Indeed, the exceptional sets
\[
\left\{\boldsymbol{\alpha}: \bigl\|\sum_{j=1}^d h_j\alpha_j\bigr\|\le |\boldsymbol h|^{-d}(\log |\boldsymbol h|)^{-1-\varepsilon}\right\}
\]
have summable measures with $\sum_{\boldsymbol h\in\mathbb Z^d} |\boldsymbol h|^{-d}(\log |\boldsymbol h|)^{-1-\varepsilon}<\infty$. Thus two distinct values of $\sum_j k_j\alpha_j$ with $\boldsymbol k$ in the enlarged box $[-2,N+2]^d$ are separated on the circle by at least the right-hand side of \eqref{prop2-gap-estimate}. Consequently any circular interval of length $O(N^{-d})$ contains at most $O_{\boldsymbol{\alpha},\varepsilon}((\log N)^{1+\varepsilon})$ such values. The symmetric difference between $[0,x)$ and $[u_1,x+u_1)$ is a union of $O(1)$ intervals of length $O(N^{-d})$, hence its total contribution is
\[
D_1\ll_{\boldsymbol{\alpha},\varepsilon}(\log N)^{1+\varepsilon}\ll(\log N)^{d-1+\varepsilon}
\]
for $d\ge2$.

It remains to control $D_2$ due to the change of the summation range. The symmetric difference between the integer box $\{1,\dots,N\}^d$ and the perturbed box $\{k_j:u_{j+1}<k_j\le N+u_{j+1}\}$ is a finite union of faces on which at least one coordinate is fixed in the set $\{-1,0, N+1, N+2\}$ and the remaining coordinates vary over intervals of length $N+O(1)$. 
After adding and subtracting the corresponding expected terms, whose signed total cancels because both discrepancies contain the same volume term $N^d x$, each face produces a discrepancy for a system of linear forms of dimension $m\le d-1$, with a fixed translation of the target interval, which can be estimated by Beck's theorem \eqref{Thm of Beck}. See the sample calculation below in detail. 

By the induction hypothesis if $m\ge2$, and by the one-dimensional Khintchine theorem if $m=1$, every such face contributes
\[
O_{\boldsymbol{\alpha},\varepsilon}((\log N)^{m+\varepsilon})\le O_{\boldsymbol{\alpha},\varepsilon}((\log N)^{d-1+\varepsilon}).
\]
There are only $O_d(1)$ faces, so the same bound holds for the whole range-perturbation error $D_2$. 

Combining the estimates for $D_1$ and $D_2$ proves \eqref{Difference between the smooth and original discrepancy}. Taking the finite intersection of the full-measure sets required in all lower dimensions gives a single full-measure set on which the estimate holds for all $N$ and all $x$.

\end{proof}

\paragraph{Sample calculation for $d=2$, $u_2=-1$, $u_3=1$.}\leavevmode\\

In this case,
\[
D(\boldsymbol{\alpha}; u_1,x+u_1; (-1,1), (N-1,N+1))
=\sum_{\substack{0\le k_1\le N-1\\ 2\le k_2\le N+1}}
\chi_{[u_1,x+u_1)}\!\left(k_1 \alpha_1 + k_2 \alpha_2 \bmod 1\right)
- N^2 x,
\]
whereas
\[
D(\boldsymbol{\alpha},x;N)
=\sum_{\substack{1\le k_1\le N\\ 1\le k_2\le N}}
\chi_{[0,x)}\!\left(k_1 \alpha_1 + k_2 \alpha_2 \bmod 1\right)
- N^2 x.
\]
Therefore,
\begin{align*}
&D(\boldsymbol{\alpha}; u_1,x+u_1; (-1,1),(N-1,N+1))
 - D(\boldsymbol{\alpha},x;N) \\
&= \sum_{\substack{0\le k_1\le N-1\\ 2\le k_2\le N+1}}
\Bigl(\chi_{[u_1,x+u_1)} - \chi_{[0,x)}\Bigr)
\!\left(k_1 \alpha_1 + k_2 \alpha_2 \bmod 1\right) \\
&\quad + \Biggl(
\sum_{2\le k_2\le N+1}
\chi_{[0,x)}\!\left(k_2 \alpha_2 \bmod 1\right)
+ \sum_{0\le k_1\le N-1}
\chi_{[0,x)}\!\left(k_1 \alpha_1 + (N+1)\alpha_2 \bmod 1\right) \\
&\qquad
- \sum_{1\le k_2\le N}
\chi_{[0,x)}\!\left(N\alpha_1 + k_2 \alpha_2 \bmod 1\right)
- \sum_{1\le k_1\le N}
\chi_{[0,x)}\!\left(k_1 \alpha_1 + \alpha_2 \bmod 1\right)
\Biggr) \\
&=: D_1 + D_2.
\end{align*}
The first term is precisely $D_1$, while the expression in parentheses is $D_2$.

By Khintchine's theorem, i.e. the case $d=1$ of \eqref{Thm of Beck}, each one-dimensional sum appearing in $D_2$ differs from its expected value by at most $O\bigl((\log N)^{1+\varepsilon}\bigr)$. For instance,
\[
\left|\sum_{2\le k_2\le N+1}
\chi_{[0,x)}\!\left(k_2 \alpha_2 \bmod 1\right)-Nx\right|
\ll_{\bma,\varepsilon} (\log N)^{1+\varepsilon}.
\]
Combining the four boundary sums gives the required estimate for $D_2$.

\section{Estimating the tail of the discrepancy function}\label{section 2 estimation of the tails}

In this section we estimate the tail of the series $\bar{D}(\bma,x;N)$, namely the contribution of Fourier modes with either very large or very small denominators. 


The reduction in this section proceeds in four successive truncations. First, Proposition \ref{Proj two D_1 - D convergence of the poisson summation} removes very large $|n_1|$ and establishes absolute convergence of the smoothed Fourier series. Second, Proposition \ref{control for sum when one coordinate is large} removes terms for which one of the coordinates $|n_1\alpha_j-n_{j+1}|$ is not the nearest-integer representative. Third, Proposition \ref{Proj two control for sum of small divisors} controls the highly resonant small-divisor part by a double-exponential decomposition. Finally, Proposition \ref{Control for sum outside Nd/4} restricts to $|n_1|\le N^d/4$, which places the main smoothing factor in a stable range for the cancellation arguments of Sections \ref{Section 3 of cancellation of the main terms.} and \ref{section 5 estimation of small exponentials}.

Recall that $\bar{D}(\bma,x;N)$ is the sum of the following terms over $\mbf{n}\in \Z^{d+1}\setminus\{\mbf{0}\}$:

\beq{\label{Proj two term}}
\begin{aligned}f(\mbf{n}, x,\bma)=
&\i^{d+1}  \frac{1-e^{2\pi \i n_1 x }}{2\pi n_1}\(\fc{\sin2\pi \(\fc{n_1}{N^d}\)}{2\pi \(\fc{n_1}{N^d}\)}\)^2.\\
&\cdot \productinDProjtwobar
\end{aligned}
\eeq
Define
$$
\bar{D}_4(\bma,x;N)=\sum_{\mbf{n}\in U_4(\bma;N)}    f(\mbf{n},x,\bma),
$$
where 
\beq\label{Proj two sum of small divisors with n1 smaller than N^d/4}
U_4(\bma;N)=
\l\{\mbf{n}\in\Z^{d+1} \ \l| \ 
\bal 
&1\le|n_1|\le N^d/4, \\
&|n_1|\prod_{j=1}^{d}\|n_1\a_j\|>(\log N)^{3d+6},\\
&|n_1\a_j-n_{j+1}|=\|n_1\a_j\|, \ 1\le j\le d.
\eal
\r.\r\}.
\eeq

Recall that $\|\cdot\|$ denotes the distance to the nearest integer. The exponent ${3d+6}$ is chosen in accordance with Lemma \ref{Proj two d_N error estimation for the number of elements}.

The main result of this section is the following:
\begin{proposition}\label{Proj two D4}
If $\varphi(n)$ is a positive nondecreasing function such that $\sum_{n=1}^{\infty} \fc{1}{\varphi(n)}<\infty$, then for almost every $\bma\in [0,1]^d$, we have 
$$|\bar{D}_4(\bma,x;N)-\bar{D}(\bma,x;N)|\ll_{\bma, \varphi} (\log N)^d\varphi^{3}(\log\log N)$$
\end{proposition}
To prove Proposition \ref{Proj two D4}, we estimate the difference $\bar{D}_4-\bar{D}$ in four steps.

\subsection{Estimating the sum when \texorpdfstring{$|n_1|$}{} is large.}

We first show that the contribution from terms with large $|n_1|$ is negligible. In particular, this proves the absolute convergence of the series \eqref{Final expression of the smoothed discrepancy} and justifies the use of Poisson's summation formula.
Define
\beq
\bar{D}_1(\bma,x;N)=\sum_{\mbf{n}\in U_1(N)}   f(\mbf{n},x,\bma),
\eeq
where 
$$U_1(N)=\{\mbf{n}\in \Z^{d+1} \backslash \{\mbf{0}\}:\l|n_1\r|<N^d(\log N)^d\}$$
We show the following:
\begin{proposition}\label{Proj two D_1 - D convergence of the poisson summation}
For almost every $\bma\in \R^d$ we have
$$\l|\bar{D}_1(\bma,x;N)-\bar{D}(\bma,x;N)\r|=\mathcal{O}(1),$$
where $\mathcal{O}(1)$ represents an absolute bound which may depend on $\bma$, but does not depend on $x$ or $N$.
\end{proposition}
\begin{proof}
In this proof we suppress subscripts in the Vinogradov symbols, since all implied constants depend at most on $\bma$.

We estimate each factor in $f(\mbf{n},x,\bma)$ according to its size.

When $|n_1\alpha_j-n_{j+1}|<1/2$ or $|n_1/N^d|<1$, that is, when the corresponding factor is \emph{small}, we use the bounds
$$\l|\fc{\sin(2\pi(n_1\alpha_j-n_{j+1}))}{2\pi(n_1\alpha_j-n_{j+1})}\r|\le 1.$$ 
and 
$$\l|\fc{\sin(2\pi(\fc{n_1}{N^d})}{2\pi(\fc{n_1}{N^d})}\r|\le 1.$$ 

When $|n_1\alpha_j-n_{j+1}|\ge 1/2$ or $|n_1/N^d|>1$, that is, when the corresponding factor is \emph{large}, we bound the numerator by $1$ and obtain
$$\l|\fc{\sin(2\pi(n_1\alpha_j-n_{j+1}))}{2\pi(n_1\alpha_j-n_{j+1})}\r|\le \fc{1}{{2\pi|n_1\alpha_j-n_{j+1}|}}.$$
and
$$\l|\fc{\sin(2\pi(\fc{n_1}{N^d}))}{2\pi(\fc{n_1}{N^d})}\r|\le \fc{1}{{2\pi|\fc{n_1}{N^d}|}}.$$

Using these inequalities, we reduce the problem to partial sums in which at most one factor $|n_1\alpha_j-n_{j+1}|$ is large.

First, we have the following upper bound; see the explanation after the definition of $U_{1,j}$:

\beq\label{Proj two D1=D11+D22}
\l|\bar{D}_1(\bma,x;N)-\bar{D}(\bma,x;N)\r|\ll\bar{D}_{1,0} +\sum_{j=1}^{d}\bar{D}_{1,j},
\eeq
where 
$$
\bar{D}_{1,0}=\sum_{\mbf{n}\in U_{1,0}(\bma;N)}\fc{N^{2d}}{|n_1|^{3}}\prod_{j=1}^{d}\fc{1}{|n_1\a_j-n_{j+1}|},
$$
where 
$$U_{1,0}(\bma;N)=\l\{\mbf{n} \in \Z^{d+1}: |n_1|>N^d(\log N)^d, \quad |n_1\a_j-n_{j+1}|<\fc{1}{2},\quad 1\le j\le d\r\},$$
and 
$$
\bar{D}_{1,j}=\sum_{\mbf{n}\in U_{1,j}(\a;N)}\fc{N^{2d}}{|n_1|^{3}}\prod_{\substack{1\le i\le d\\ i\neq j}}\fc{1}{|n_1\a_i-n_{i+1}|}\cdot \fc{1}{|n_1\a_j-n_{j+1}|^{3}},
$$
where
$$
U_{1,j}({\bm{\a}};N)=\l\{\mbf{n} \in \Z^{d+1}: |n_1|>N^d(\log N)^d, \quad |n_1\a_j-n_{j+1}|>\fc{1}{2}, \quad |n_1\a_i-n_{i+1}|<\fc{1}{2}, \ i\neq j\r\}.
$$
The bound \eqref{Proj two D1=D11+D22} omits the partial sums corresponding to terms with \emph{two or more} large factors $|n_1\alpha_j-n_{j+1}|$. Indeed, 
every additional large coordinate contributes a summable factor $\sum_{r\ge1}r^{-3}$, and hence those partial sums are dominated by the case with only one large coordinate, as in \eqref{inequality to ignore n such that two or more factors are big} below.
Define $r_j=\lceil |n_1\a_j-n_{j+1}|\rceil$ for $1\le j\le d$. Then for each $1\le j\le d$ we have
\beq\label{inequality to ignore n such that two or more factors are big}
\bar{D}_{1,j}\ll\sum_{r=1}^{\infty}\sum_{\mbox{\scriptsize$\ba{c} \mbf{n} \in \Z^{d+1}:\\|n_1|>N^d(\log N)^d \\|n_1\a_i-n_{i+1}|< 1/2\\ i\neq j \ea$}}\fc{N^{2d}}{|n_1|^{3}}\prod_{\substack{1\le i\le d\\ i\neq j}}\fc{1}{|n_1\a_i-n_{i+1}|}\fc{1}{r^{3}}\ll\bar{D}_{1,0},
\eeq

Therefore, it remains to prove that $\bar{D}_{1,0}=\mathcal{O}(1)$. For this, we use the following lemma, whose proof is a standard application of the Borel--Cantelli lemma; for completeness, we include the details in Appendix \ref{appendix:proof-of-proj-two-lem1}.
\begin{lemma}\label{Proj two lem1}
If $\varphi(n)$ is a positive nondecreasing function such that $\sum_{n=1}^{\infty} \fc{1}{\varphi(n)}<\infty$, then for almost every $\bma\in \R^d $, the sum
$$\sum_{n=2}^\infty \fc{1}{n\varphi(\log n)\prod_{j=1}^{d} \l(\|n \a_j\|\varphi(|\log \|n\a_j\||)\r)}$$
converges.
\end{lemma}
Continuing with the proof for Proposition \ref{Proj two D_1 - D convergence of the poisson summation}, since $\sums_{n=1}^{\infty} {1}/{n^2}=\mathcal{O}(1)$, by the Borel--Cantelli lemma, we have for almost every $\bma \in \R^d$,
$$\|n\a_j\|\ge\fc{1}{n^2}, \quad 1\le j\le d,$$ 
for all but finitely many $n\in \N$, or equivalently,
\beq\label{Proj two inequality for log nalpha and log n}
\l|\log\|n\a_j\|\r|\le 2\log n, \quad 1\le j\le d, 
\eeq
for all but finitely many $n\in \N^+$.
Taking the function in Lemma \ref{Proj two lem1} above to be $\varphi_0(n)=n^{1+\e}$ where $\e>0$ is small, we estimate $\bar{D}_{1,0}$ as follows:
\begin{equation}\label{Proj two D10 estimate by Lemma 5}
\begin{aligned}
\bar{D}_{1,0}
&=\sum_{n_1>N^d(\log N)^d} \fc{1}{(\fc{n_1}{N^d})^2} \fc{1}{n_1\prod_{j=1}^{d} \|n_1\a_j\|}\\
&\ll\sum_{n\ge N^d(\log N)^d} \fc{1}{(\log n)^{(d+1)(1+\e)}}\fc{1}{n\prod_{j=1}^{d} \|n\a_j\|}\qquad \text{using range of } n_1\\
&\ll\sum_{n\ge N^d(\log N)^d} \fc{1}{|n|\varphi_0(\log n)\prod_{j=1}^{d} \l(\|n \a_j\|\varphi_0(|\log \|n\a_j\||)\r)}\qquad \text{using \eqref{Proj two inequality for log nalpha and log n}}\\
&\le\sum_{n=2}^{\infty}\fc{1}{|n|\varphi_0(\log n)\prod_{j=1}^{d} \l(\|n \a_j\|\varphi_0(|\log \|n\a_j\||)\r)}\\
&=\mathcal{O}(1),
\end{aligned}
\end{equation}

finishing the proof.
\end{proof}

\subsection{Estimating the sum when \texorpdfstring{$|n_1 \a_j-n_{j+1}|$}{}  \texorpdfstring{$> 1/2$}{} for some \texorpdfstring{$1\le j\le d$}{}.}
Let 
$$
\bar{D}_2(\bma,x;N)=\sum_{\mbf{n}\in U_2(\bma;N)}  f(\mbf{n},x,\bma)
$$
where $\mbf{n}=(n_1,\dots, n_{d+1})$ satisfies
\beq{\label{Proj two bigger than one third}}
U_2(\bma;N)=
\l\{\mbf{n}\in\Z^{d+1}\backslash \{\mbf{0}\} \ \l| \ 
\bal 
&|n_1|<N^d(\log N)^d,\\
&|n_1\a_j-n_{j+1}|=\|n_1\a_j\|, \ 1\le j\le d\\
\eal
\r.\r\}
\eeq
We will show that $\bar{D}_1$ can be replaced by $\bar{D}_2$:
\begin{proposition}\label{control for sum when one coordinate is large}
For almost every $\bma\in [0,1]^d$, we have
$$
\l|\bar{D}_2(\bma,x;N)-\bar{D}_1(\bma,x;N)\r| \ll_{\bma} (\log N)^d
$$
\end{proposition}

\begin{proof}
Since we can bound those $\mbf{n}\in U_1(\bma;N)$ for which $|n_1\alpha_j-n_{j+1}|\ge {1}/{2}$ holds for two or more indices by the corresponding sums where only one coordinate is treated as large, as in \eqref{inequality to ignore n such that two or more factors are big}, we have
\beq\label{Proj two sum for small divisor n||n*alpha||}
\bal
\l|\bar{D}_2(\bma,x;N)-\bar{D}_1(\bma,x;N)\r|\ll 
&\sum_{j=1}^{d}\sum_{\mbf{n}\in U_{2,j}}\fc{1}{|n_1|}\prod_{\substack{1\le i\le d\\ i\neq j}}\fc{1}{|n_1\a_i-n_{i+1}|}\cdot \fc{1}{|n_1\a_j-n_{j+1}|^{3}}\\
&\ll \sum_{j=1}^{d}\sum_{n=1}^{N^d\l(\log N\r)^d} \fc{1}{n\prod_{i\neq j}\|n\a_i\|},
\eal
\eeq
where 
$$
U_{2,j}=\l\{|n_1|\le N^d(\log N)^d, \quad |n_1\a_i-n_{i+1}|< \fc{1}{2}, \ i\neq j, \quad |n_1\a_j-n_{j+1}|\ge \fc{1}{2}\r\}.
$$

We use the following direct consequence of \cite[Theorem~1.4]{fregoli2024sumsreciprocalsfractionalparts}. In Fregoli's notation this is the almost-sure bound for the averaged reciprocal sum $S^*(\boldsymbol\beta,T)$; taking the symmetric parameter $T=N^{d}(\log N)^{d}$ gives exactly the form below.

\begin{lemma}\label{fregoli-reciprocal-input}
For almost every $\bm{\beta}=(\beta_1,\dots, \beta_{d-1} )\in \R^{d-1}$ and every $N>1$, we have 
$$
\sum_{n=1}^{N^{d}(\log N)^{d}} \fc{1}{n \prod_{i=1}^{d-1}\|n\beta_i\|} \ll_{\bm{\beta}} (\log N)^d.
$$
\end{lemma}
Since the sum in \eqref{Proj two sum for small divisor n||n*alpha||} is a finite sum over $1\le j\le d$, the proof of the proposition is completed.

\end{proof}
\begin{remark}
As an alternative to \cite{fregoli2024sumsreciprocalsfractionalparts}, Appendix \ref{appendix:weak-fregoli-substitute} gives a direct short proof of a slightly weaker bound, which in the present application may be written as $(\log N)^d\varphi^{d+1}(\log\log N)$ after enlarging the secondary logarithmic exponent. This keeps the optimal main term and only weakens the second-order factor. The proof uses Lemma \ref{Proj two lem1}, with an auxiliary function of the form $\varphi_0(t)=t\varrho(\log t)$ as in \eqref{Proj two D10 estimate by Lemma 5}, and the same replacement made for the factors involving $|\log\|n\beta_i\||$.
\end{remark}

\subsection{Estimating the sum when \texorpdfstring{$|n_1|\prod_{j=1}^{d} \|n_1\a_j\|$}{} is small}

Define
\beq\label{Proj two sum of small divisors}
\bar{D}_3(\bma,x;N)=\sum_{\mbf{n}\in U_3(\bma;N)}  f(\mbf{n},x,\bma),
\eeq
where
\beq\label{Proj two U3}
U_3(\bma;N)=
\left\{\mbf{n}\in\Z^{d+1}\setminus \{\mbf{0}\} \ \middle| \ 
\begin{aligned}
&1\le |n_1| \le N^d(\log N)^d,\\
&|n_1| \prod_{j=1}^{d} \|n_1\a_j\|> (\log N)^{3d+6},\\
&|n_1\a_j-n_{j+1}|=\|n_1\a_j\|, \quad 1\le j\le d
\end{aligned}
\right\}.
\eeq

Observe that in \eqref{Proj two U3}, each $n_{j+1}$ is the closest integer to $n_1\a_j$. Hence, $\mbf{n}$ is uniquely determined by $n_1$ and $\bma$. In the sequel, it therefore suffices to consider $n_1$ instead of $\mbf{n}\in \Z^{d+1}$.

The main result of this step is the following proposition.

\begin{proposition}\label{Proj two control for sum of small divisors}
\convergentconditionforphi. Then, for almost every $\bma$, we have
$$
\left|\bar{D}_3(\bma,x;N)-\bar{D}_2(\bma,x;N) \right|\ll_{\bma, \varphi} (\log N)^d\varphi^{3}(\log \log N).
$$
\end{proposition}
\begin{proof}
Observe that
$$
\left|\bar{D}_3(\bma,x;N)-\bar{D}_2(\bma,x;N)\right|
\ll
\sum_{\substack{
1\le n \le N^d(\log N)^d \\
n\prod_{j=1}^{d} \|n\a_j\|<(\log N)^{3d+6}
}}
\frac{1}{n\prod_{j=1}^{d} \|n\a_j\|}.
$$
So the goal is to show that:
\beq\label{Proj two lemma for sum of small divisors}
\sum_{\substack{
1\le n \le N^d(\log N)^d \\
n\prod_{j=1}^{d} \|n\a_j\|<(\log N)^{3d+6}
}}
\frac{1}{n\prod_{j=1}^{d} \|n\a_j\|}
\ll_{\bma, \varphi} (\log N)^d\varphi^{3}(\log \log N).
\eeq
We first apply two minor modifications to our function $\varphi$.
Fix a sufficiently large constant $A=A(d)\ge 1$ and put, for $t\ge 1$,
$$
        \widetilde{\varphi}(t)=\varphi\big(\max\{\lfloor t/A\rfloor,1\}\big),\qquad
        \psi(t)=\min\{\widetilde{\varphi}(t),t^2\}.
$$
Then $\psi$ is positive and nondecreasing, $\psi(t)\le \varphi(t)$, $\psi(t)\le t^2$, and
$$
        \sum_{m=1}^{\infty}\frac{1}{\psi(m)}<\infty.
$$
The harmless dilation by $A$ is introduced only to ensure that, whenever $p$ or $v$ is bounded by a dimensional multiple of $\log\log N$, the corresponding value of $\psi$ is bounded by $\varphi(\log\log N)$. This is used in \eqref{usage of dilation}. The minimum of $\{\widetilde{\varphi}(t),t^2\}$ is to control the range of $v$ to be of order $\log\log N$ if $\varphi$ grows too rapidly.

We partition the integers $n$ into the sets
$$
S_{\bma}(p,v)=\left\{
 n\in \mathbb{N}: e^{e^{p-1}}\le n<e^{e^p},\quad
 \frac{2^{v-1}}{(\log n)^{d}\psi(\log \log n)}
 \le n\prod_{j=1}^{d} \|n\a_j\|
 <\frac{2^v}{(\log n)^{d}\psi(\log \log n)}
\right\},
$$
where $\log\log n$ is interpreted as $1$ for $n<e^e$.

We first record the lower bound for the divisor, which follows from the standard cusp estimate \eqref{standard cusp estimate} and an application of the convergent Borel-Cantelli Lemma. It will be used only after the uniform cardinality estimate below has already been proved.

\begin{lemma}\label{Proj two absolute lower bound the the divisor}
For a positive nondecreasing function $\psi(n)$ such that $\sum_{n=1}^{\infty}1/\psi(n)<\infty$, for almost every $\bma\in [0,1]^d$, and for all sufficiently large $n$, we have
$$
        n\prod_{j=1}^{d} \|n\a_j\|\gg_{\bma,\psi}\frac{1}{(\log n)^{d}\psi(\log \log n)},
$$
and, for each $1\le j\le d$,
$$
        \|n\a_j\|\gg_{\bma,\psi} \frac{1}{n\psi(\log n)}.
$$
\end{lemma}

The next proposition is proved for all $v\in\mathbb{Z}$ independent of $\bma$.

\begin{proposition}\label{cardinality for Spv}
Let $v_-=\max\{-v,0\}$ and $\langle v\rangle=\max\{|v|,1\}$. For almost every $\bma$, for all $p>1$ and all $v\in\mathbb{Z}$, 
one has
\beq\label{cardinality for Spv corrected}
        \#S_{\bma}(p,v)
        \ll_{\bma,\varphi,d}
        2^v\psi(\langle v\rangle)(1+v_-)^{d-1}.
\eeq
\end{proposition}

\begin{proof}
Define
$$
C(n,v)=\left\{\bma\in[0,1]^d:
\frac{2^{v-1}}{(\log n)^d\psi(\log\log n)}
\le n\prod_{j=1}^{d}\|n\a_j\|
<\frac{2^v}{(\log n)^d\psi(\log\log n)}
\right\}.
$$
Put
$$
        \tau_{n,v}=\frac{2^v}{n(\log n)^d\psi(\log\log n)}.
$$
The standard multiplicative-cusp estimate gives, uniformly in $0<t$,
\beq\label{standard cusp estimate}
        \operatorname{Leb}\left\{y\in[0,1]^d:\prod_{j=1}^{d}\|y_j\|<t\right\}
        \ll_d t\bigl(1+\log^+(1/t)\bigr)^{d-1},
\eeq
where $\log^+ u=\max\{\log u,0\}$. Since $n\bma$ is measure preserving modulo one,
\beq\label{cusp estimate with negative v factor}
        \operatorname{Leb}(C(n,v))\ll_d
        \frac{2^v}{n(\log n)^d\psi(\log\log n)}
        \bigl(\log n+1+v_-\bigr)^{d-1}.
\eeq
Here we used $\psi(t)\le t^2$, so that $\log\psi(\log\log n)\ll \log\log\log n\ll \log n$; for negative $v$ the only additional possible growth is precisely the factor $v_-^{d-1}$.

For $e^{e^{p-1}}\le n<e^{e^p}$, we have $\log\log n\in[p-1,p]$. Hence, for $p\ge2$,
\[
\begin{aligned}
\mathbb{E}\#S_{\bma}(p,v)
&=\sum_{e^{e^{p-1}}\le n<e^{e^p}}\operatorname{Leb}(C(n,v))  \\
&\ll_d \frac{2^v}{\psi(p-1)}
\sum_{e^{e^{p-1}}\le n<e^{e^p}}\frac{(\log n+1+v_-)^{d-1}}{n(\log n)^d}  \\
&\ll_d \frac{2^v(1+v_-)^{d-1}}{\psi(p-1)} .
\end{aligned}
\]
Indeed, after the substitution $u=\log n$, the remaining sum is bounded by
$$
        \int_{e^{p-1}}^{e^p}\frac{(u+1+v_-)^{d-1}}{u^d}\,du
        \ll_d (1+v_-)^{d-1}.
$$
By Markov's inequality,
$$
\operatorname{Leb}\left\{\#S_{\bma}(p,v)>C2^v\psi(\langle v\rangle)(1+v_-)^{d-1}\right\}
\ll_d \frac{1}{\psi(p-1)\psi(\langle v\rangle)}.
$$
Since
$$
        \sum_{p=2}^{\infty}\sum_{v\in\mathbb{Z}}\frac{1}{\psi(p-1)\psi(\langle v\rangle)}<\infty,
$$
the Borel--Cantelli lemma proves \eqref{cardinality for Spv corrected} for all but finitely many pairs $(p,v)$. Enlarging the implied constant, which is allowed to depend on $\bma$, absorbs these finitely many exceptions.
\end{proof}

We now complete the proof of Proposition \ref{Proj two control for sum of small divisors}. Fix $\bma$ in the full-measure set for both Lemma \ref{Proj two absolute lower bound the the divisor} and Proposition \ref{cardinality for Spv}. The lemma gives a constant $K=K(\bma,\varphi)$ such that for all sufficiently large $n$, only $v\ge -K$ occur; the finitely many remaining $n$ are absorbed into the implied constant.

If $n\le N^d(\log N)^d$, then $p\le \log\log N+O_d(1)$. Moreover, the additional restriction
$$
        n\prod_{j=1}^{d}\|n\a_j\|<(\log N)^{3d+6}
$$
implies, because $\psi(t)\le t^2$, that $v\le C_d\log\log N$. By choosing the constant $A=A(d)$ above sufficiently large, all positive values of $p$ and $v$ which occur in the following sums satisfy
\begin{equation}\label{usage of dilation}
       \psi(p)\le \varphi(\log\log N),\qquad
        \psi(\langle v\rangle)\le \varphi(\log\log N).
\end{equation}
   
Therefore, using Proposition \ref{cardinality for Spv},
\beq
\begin{aligned}
&\sum_{\substack{
1\le n \le N^d(\log N)^d \\
n\prod_{j=1}^{d} \|n\a_j\|<(\log N)^{3d+6}
}}
\frac{1}{n\prod_{j=1}^{d} \|n\a_j\|}\\
&\ll_{\bma,\varphi}
\sum_{p\le \log\log N+O_d(1)}
\sum_{-K\le v\le C_d\log\log N}
\frac{e^{dp}\psi(p)}{2^v}
\#S_{\bma}(p,v)\\
&\ll_{\bma,\varphi}
\sum_{p\le \log\log N+O_d(1)} e^{dp}\psi(p)
\sum_{-K\le v\le C_d\log\log N}\psi(\langle v\rangle)(1+v_-)^{d-1}\\
&\ll_{\bma,\varphi}
(\log N)^d\varphi(\log\log N)\bigl(\varphi(\log\log N)\log\log N\bigr)\\
&\ll_{\bma,\varphi}(\log N)^d\varphi^3(\log\log N).
\end{aligned}
\eeq
In the last line we used the elementary consequence of $\sum_m1/\varphi(m)<\infty$ and monotonicity, namely $m\ll_{\varphi}\varphi(m)$ for all sufficiently large $m$. This proves \eqref{Proj two lemma for sum of small divisors}, and thus the proof is completed.
\end{proof}

\subsection{Controlling the sum when \texorpdfstring{$n_1$}{} lies between \texorpdfstring{$N^d/4$}{} and \texorpdfstring{$N^d(\log N)^d$}{}.}
The goal of this step is to restrict the range of $2\pi n_1/N^d$ in $f(\mbf{n},x,\bma)$ to $(-\pi/2,\pi/2)$. This gives better control of $f$ and will be useful when we later estimate the cancellation among the main terms. We prove the following proposition.

Recall
\beq
\bar{D}_4(\bma,x;N)=\sum_{\mbf{n}\in U_4(\bma;N)}  f(\mbf{n},x,\bma),
\eeq
where
\beq\label{Proj two U4}
U_4(\bma;N)=
\left\{\mbf{n}\in\Z^{d+1}\setminus \{\mbf{0}\} \ \middle| \ 
\begin{aligned}
&1\le |n_1| \le \frac{N^d}{4},\\
&|n_1| \prod_{j=1}^{d} \|n_1\a_j\|> (\log N)^{3d+6},\\
&|n_1\a_j-n_{j+1}|=\|n_1\a_j\|, \quad 1\le j\le d.
\end{aligned}
\right\}.
\eeq

\begin{proposition}\label{Control for sum outside Nd/4}
For almost every $\bma $, we have 
$$\l|\bar{D}_4(\bma,x;N)-\bar{D}_3(\bma,x;N)\r|\ll_{\bma}(\log N)^d \log \log N.
$$
\end{proposition}

We decompose the range of $n_1$ appearing in $U_3\setminus U_4$ as follows. Let
\beq\label{decomposition of the set between Nd/4 to NdlnNd}
T_{\bma}(\mbf{l};N)=
\l\{
n_1 \in \Z \l|  \
\bal
&2^{l_1}\le |n_1|< 2^{l_1+1},\\
&2^{-l_{j+1}}\le \|n_1\a_j\|<2^{-l_{j+1}+1},\ 1\le j\le d-1\\
&2^{l_{d+1}-l_1+\sum_{i=1}^{d-1}l_{i+1}}\le \|n_1\a_d\|< 2^{l_{d+1}-l_1+\sum_{i=1}^{d-1}l_{i+1}+1}\\
\eal
\r.\r\},
\eeq
Here $l_1$ specifies
 the dyadic range of $|n_1|$. Since $|n_1|\le N^{d}(\log N)^{d}$ in our setting, we have
$$
2^{l_1}\le N^{d}(\log N)^{d}.
$$
For each $1\le j\le d-1$, the index $l_{j+1}$ specifies the dyadic range of $\|n_1\a_j\|$. By Lemma \ref{Proj two absolute lower bound the the divisor}, these quantities cannot be too small, and hence $l_{j+1}\ll \log N$. Finally, $l_{d+1}$ specifies the dyadic range of
$$
|n_1|\prod_{j=1}^{d}\|n_1\a_j\|,
$$
and from the defining condition of this case we obtain
$$
2^{l_{d+1}}\ge (\log N)^{3d+6}.
$$

The admissible set of indices $\mbf{l}=(l_1,\dots,l_{d+1})\in L_1(N)$ is given by the following bounds. 
\beq\label{Proj two range for l between N^d/4 and N^d (log N)^2}
L_1(N):=\l\{\mbf{l}\in \Z^{d+1} \l| \
\bal
&d\log_2 N-2\le l_1\le d\log_2 N +d\log_2\log N+2, \\
&2 \le l_{j+1}\ll_{\bma} \log N, \ 1\le j\le d-1,\\ 
&(\log N)^{3d+6} \le 2^{l_{d+1}}\le 2^{l_1}\le N^{d+1},
\eal
\r. \r\}.
\eeq
We first prove
 an upper bound with a uniform constant for $\#T_{\bma}(\mbf{l};N)$.
\begin{lemma}\label{Proj two number of elements in S for N^d/4}
For almost every $\bma\in \R^d$ and all $\mbf{l}\in L_1(N)$, we have
$$\#T_{\bma}(\mbf{l};N)\ll_{\bma} 2^{l_{d+1}}$$
where the constant is uniform for $\mbf{l}$ in $L_1(N)$.
\end{lemma}
\begin{proof}
The cardinality of $T_{\bma}(\mbf{l};N)$ is the number of integers $n_1\in[2^{l_1},2^{l_1+1})$ such that the point
$$n_1\bma \mod 1=(n_1\a_1 \mod 1,\dots,n_1\a_d\mod 1)$$
lies in one of the $2^d$ target boxes, and we denote their union by $B(\mbf{l})$,
$$\prod_{j=2}^{d}\pm[2^{-l_{j}},2^{-l_{j}+1})\times \pm[2^{l_{d+1}-l_1+\sum_{i=1}^{d-1}l_{i+1}},2^{l_{d+1}-l_1+\sum_{i=1}^{d-1}l_{i+1}+1}),$$ 
where the negative signs account for the case $\fc{1}{2}\le \{n_1\a_j\}<1$.

We use Beck's theorem \eqref{Thm of Beck} to estimate $\#T_{\bma}(\mbf{l};N)$. First observe that
\beq\label{expression for Talpha L N}
\#T_{\bma}(\mbf{l};N)=\sum_{n=2^{l_1}}^{2^{l_1+1}-1} \chi_{B(\mbf{l})}(n\bma \mod 1)= \sum_{n=1}^{2^{l_1}} \chi_{B(\mbf{l})-(2^{l_1}-1)\bma}(n\bma \mod 1)
\eeq

The boxes $\chi_{B(\mbf{l})-(2^{l_1}-1)\bma}$ are translated on the torus,
and for fixed $\alpha$, each box can be written as a finite signed sum of boxes with corners at the origin, so Beck's estimate \eqref{Thm of Beck} applies for a full-measure set of $\a$.

Note that each of the $2^d$ boxes has the volume $2^{l_{d+1}-l_1}$, thus 
\[\text{Vol}(B(\mbf{l}))=\mathcal{O}(2^{l_{d+1}-l_1}).\]

By \eqref{expression for Talpha L N}, take in Beck's estimate \eqref{Thm of Beck} $\varphi(n)=n^{2}$ and use $(\log\log n)^2\ll (\log n)^{1/2}$, and we have that for almost every $\bma$, 
\[
\begin{aligned}
\#T_{\bma}(\mbf{l};N)
&=\sum_{n=1}^{2^{l_1}} \chi_{B(\mbf{l})-(2^{l_1}-1)\bma}(n\bma \mod 1)\\
&=2^{l_1} \text{Vol}(B(\mbf{l})) + \mathcal{O}((\log 2^{l_1})^{d+\frac{1}{2}})\\
&=\mathcal{O}({2^{l_{d+1}}})+\mathcal{O}(l_1^{d+\frac{1}{2}}).
\end{aligned}
\]

It is important to carefully check that the constant in Beck's bound is uniform for all the boxes in $B(\mbf{l})$ and all the $\mbf{l}\in L_1(N)$ (see \eqref{Proj two range for l between N^d/4 and N^d (log N)^2}). Indeed, by definition of the maximal discrepancy \eqref{Beck maximal discrepancy}, the bound is uniform for all boxes inside $B(\mbf{l})$. Furthermore, the constant in Beck's bound \eqref{Thm of Beck} is uniform in $N$, thus changing $\mbf{l}$ does not compromise the bound.


Note that $l_1^{d+\fc{1}{2}}\ll(\log N)^{d+\fc{1}{2}}<(\log N)^{3d+6}<2^{l_{d+1}}$, so the claim follows.
\end{proof}
\begin{remark}
  Here we note that it is essential that the volume of the dyadic boxes appearing in $B(\mbf{l})$ is not too small, and that the same argument will reappear in Lemma \ref{Proj two d_N error estimation for the number of elements}. The reason this argument can replace a second moment estimate is precisely the duality between Kronecker sequences and linear forms. In fact, one could exploit the duality in the other direction and use a weaker bound for linear forms, such as $(\log N)^{d+2}$, to prove a slightly weaker upper bound in Beck's theorem \eqref{Thm of Beck} without a second moment argument. 
\end{remark}

We can now estimate the contribution of the terms with $N^d/4\le |n_1|\le N^d(\log N)^d$.

\nid \textbf{Proof of Proposition \ref{Control for sum outside Nd/4}}
\begin{proof}
Suppose that $\a$ satisfies Lemma \ref{Proj two number of elements in S for N^d/4}. Inside each set $T_{\bma}(\mbf{l};N)$, the divisor $|n_1|\prod_{j=1}^{d}\|n_1\a_j\|$ lies between $2^{l_{d+1}}$ and $2^{l_{d+1}+d+1}$,
$$\bal
&\quad \l|\bar{D}_4-\bar{D}_3\r|\\
&\le \sum_{\mbf{l}\in L_1(N)}\sum_{T_{\bma}(\mbf{l};N)} \fc{1}{|n_1|\prod_{j=1}^{d} \|n_1\a_j\|}\\
&\le \sum_{\mbf{l}\in L_1(N)}\sum_{T_{\bma}(\mbf{l};N)} \fc{1}{2^{l_{d+1}}}\\
&\ll_{\bma} \sum_{\mbf{l}\in L_1(N)}\fc{1}{2^{l_{d+1}}}\#T_{\bma}(\mbf{l};N)\\
&\ll \sum_{\mbf{l}\in L_1(N)} 1\\
&\ll (\log N)^d \log \log N.
\eal$$
The last inequality follows from the definition of $L_1(N)$ in \eqref{Proj two range for l between N^d/4 and N^d (log N)^2}: there are $\mathcal{O}(\log\log N)$ possible values of $l_1$, and $\mathcal{O}(\log N)$ possible values of each $l_i$ for $2\le i\le d+1$.
\end{proof}

\noindent\textbf{Proof of Proposition \ref{Proj two D4}}
\begin{proof}
Combining Propositions \ref{Proj two D_1 - D convergence of the poisson summation}, \ref{control for sum when one coordinate is large}, \ref{Proj two control for sum of small divisors}, \ref{Control for sum outside Nd/4}, we obtain
$$
\bal 
|\bar{D}_4-\bar{D}| 
&\ll_{\bma,\varphi} \mathcal{O}(1)+(\log N)^d+(\log N)^d\varphi^{3}(\log \log N)+
(\log N)^d \log \log N\\
&\ll_{\bma,\varphi} (\log N)^d\varphi^{3}(\log \log N).\\
\eal
$$
\end{proof}

In Section \ref{Section 3 of cancellation of the main terms.} we control the non-oscillatory part by pairing signed cells whose divisors have opposite signs.  In Section \ref{section 5 estimation of small exponentials} we then control the oscillatory parts by showing that, for each fixed exponential phase, the cells without sufficient cancellation lie on sparse special lines.

\section{Cancellation of the main terms}\label{Section 3 of cancellation of the main terms.}
To estimate the contribution of the remaining "main" terms, 
we first decompose the product $f(\mbf{n},x,\bma)$ into two types of terms, according to the appearance of exponentials in the numerator:
\beq\label{Proj two multiplied out term}
\bal
\bar{D}_4=\fc{\i^{d+1}}{(2\pi)^{d+1}}
\Bigg(
&\sum_{\mbf{n}\in  U_4(\bma;N)} \fc{1}{n_1\prod_{j=1}^{d} (n_1\a_j-n_{j+1})}\cdot g(\bm{n}, \bma; N)\\
&+\sum_{\mbf{s}}\pm\sum_{\mbf{n}\in U_4(\bma;N)} \fc{e^{2\pi \i \L_{\mbf{s}}(\mbf{n})}}{n_1\prod_{j=1}^{d} (n_1\a_j-n_{j+1})}\cdot g(\bm{n}, \bma; N)\Bigg).\\
\eal
\eeq
Recall $U_4$ from \eqref{Proj two U4}. Here $g(\mbf{n}, \bma; N)$ is the product below (observe that $|g(\mbf{n},\bma; N)|\le 1$):
\beq\label{Proj two gN}
\l(\fc{\sin 2\pi(\fc{n_1}{N^d})}{2\pi(\fc{n_1}{N^d})}\r)^2\prod_{j=1}^{d}\l(\fc{\sin 2\pi (n_1 \a_j -n_{j+1})}{2\pi (n_1\a_j-n_{j+1})}\r)^2
\eeq
and $\L_{\mbf{s}}=\L_{\mbf{s},x,N,\bma}$ is one of the $2^{d+1}-1$ linear forms of $d+1$ variables:
\beq\label{Proj two linear forms}
\L_{\mbf{s}}(\mbf{n})=\L_{\mbf{s}}(n_1,\dots,n_{d+1})=\,\sigma_1 n_1 x-
\sum_{j=1}^{d}\sigma_{j+1}N(n_1\a_j-n_{j+1}).
\eeq
The index $\mbf{s}=(\sigma_1,\dots,\sigma_{d+1})\in\{0,1\}^{d+1}$ and $\mbf{s}\neq\mbf{0}$.
The sign $\pm$ in the second part of \eqref{Proj two multiplied out term} is
$\pm=(-1)^{\sum_{j=1}^{d+1}\sigma_j}$, and hence it is independent of $\mbf{n}\in\Z^{d+1}$.
In this section, we control the constant part by pairing signed cells $S_{\bma}(\mbf{l},\signvec^{\pm};N)$ (defined below), whose divisors have opposite signs. More precisely, the dyadic cells are paired so that their expected cardinalities match, whereas the signs of the divisors are opposite; the resulting cancellation gives the sufficient logarithmic saving.  


We begin with the (constant) part of $\bar{D}_4$.

Let
\beq\label{Proj two sum9}
\bar{D}_{5}(\bma,x;N)=\sum_{\mbf{n}\in U_4(\bma;N)} \fc{g(\bm{n}, \bma; N)}{n_1\prod_{j=1}^{d} (n_1\a_j-n_{j+1})}.
\eeq
The goal of this section is to prove the following:
\begin{proposition}\label{Proj two Prop cancellation of the main terms}
For almost every $\bma\in [0,1]^d$, we have 
$$
\l|\bar{D}_{5}(\bma,x;N)\r|\ll_{\bma} (\log N)^{d-1}
$$
\end{proposition}

The idea is to decompose the Fourier modes into sets where all the factors of the divisor are basically of the same magnitude, and only one factor differs in sign ($+$ or $-$). 
Pairing the positive cells with the corresponding negative cells yields the cancellation that gives the saving of $\log N$.

Let $\d_N=1/(\log N)^2$, $\signvec=(\epsilon_1,\dots,\epsilon_{d+1})\in\{\pm1\}^{d+1}$, and $\mbf{l}=(l_1,\dots,l_{d+1})\in \Z^{d+1}$. We define the following cells, in analogy with \eqref{decomposition of the set between Nd/4 to NdlnNd}:
\beq\label{Proj two def for S}
S_{\bma}(\mbf{l},\signvec;N)=
\l\{
\bal
&\mbf{n} \in \Z^{d+1}: (1+\d_N)^{l_1}\le \epsilon_1 n_1< (1+\d_N)^{l_1+1},\\
&(1+\d_N)^{-l_{j+1}}\le\epsilon_{j+1} (n_1 \a_j-n_{j+1})<(1+\d_N)^{-l_{j+1}+1},\ 1\le j\le d-1,\\
&(1+\d_N)^{l_{d+1}-l_1+\sum_{i=1}^{d-1}l_{i+1}}\le \epsilon_{d+1} (n_1\a_{d}-n_{d+1})\le (1+\d_N)^{l_{d+1}-l_1+\sum_{i=1}^{d-1}l_{i+1}+1}, \\
& \text{each $n_{j+1}$ is the closest integer to $n_1 \a_j$, $1\le j\le d$}.\\
\eal
\r\},
\eeq
where the $l_i$ are positive integers, and the range of $\mbf{l}=(l_1, \dots, l_{d+1})$ is as the following:
$$
L_2(N)=\l\{ \mbf{l}\in \N^{d+1} \l|\
\bal
&(\log N)^{3d+6}\le (1+\d_N)^{l_{d+1}}\le (1+\d_N)^{l_1}\le N^d/4,\\
& (1+\d_N)^{l_1-l_{j+1}}\ge(\log N)^{3d+6}, \ 1\le j\le d-1.
\eal
\r.
\r\}
$$
which gives
\beq\label{Proj two condition for l}
\log \log N/ \delta_N\ll l_{d+1} \le l_1\ll \log N/\delta_N, \quad  l_1-l_{j+1}\gg \log \log N/\d_N, \ 1\le j\le d-1.
\eeq

For each $\mbf l$ we fix \emph{once and for all} the following partition $\mathcal P$ of $\{\pm1\}^{d+1}$ into $2^d$ unordered pairs: a sign vector is paired with the vector obtained by flipping its first coordinate. Since the sign of the divisor $n_1\prod_{j=1}^d(n_1\alpha_j-n_{j+1})$ is $\prod_{i=1}^{d+1}\epsilon_i$ on $S_{\bma}(\mbf l,\signvec;N)$, the two cells in each pair have opposite signs. We write a pair as $\mathfrak p=(\signvec^+,\signvec^-)$, ordered so that the divisor is positive on $S(\mbf l,\signvec^+;N)$ and negative on $S(\mbf l,\signvec^-;N)$. All sums over paired sign choices below refer to this fixed partition and hence do not overcount sign cells.

By integration over $\bma\in [0,1)^d$, the expected value for $\#S_{\bma}(\mbf{l},\signvec;N)$ is:
$$
E(\mbf{l},\signvec;N)=\int_{\bma\in [0,1]^d} \#S_{\bma}(\mbf{l},\signvec;N)\, d\bma. 
$$
A more precise description for the number of elements in $\#S_{\bma}(\mbf{l},\signvec;N)$ is the following:
\begin{lemma}\label{Proj two d_N error estimation for the number of elements}
For almost every $\bma\in \mathbb R^d$, and every $\mbf{l}\in L_2(N)$, we have
\[
|\#S_{\bma}(\mbf{l},\signvec;N)-E(\mbf{l},\signvec;N)|\ll_{\bma} \d_N E(\mbf{l},\signvec;N),
\]
where the constant is uniform for $\mbf{l}\in L_2(N)$, $N$, and the sign vector $\signvec$.
Moreover, for every paired sign choice $\mathfrak p=(\signvec^+,\signvec^-)\in\mathcal P$, we have
\[
E(\mbf{l},\signvec^+;N)=E(\mbf{l},\signvec^-;N)\asymp_d \d_N^{d+1}(1+\d_N)^{l_{d+1}}.
\]
\end{lemma}
\begin{proof}
The set $S_{\bma}(\mbf l,\signvec;N)$ is counted by those integers
\[
(1+\d_N)^{l_1}\le \epsilon_1n_1<(1+\d_N)^{l_1+1}
\]
for which $n_1\bma$ lies in a union of $O_d(1)$ boxes on the torus with side lengths
\[
\d_N(1+\d_N)^{-l_{2}},\ldots,\d_N(1+\d_N)^{-l_d},
\d_N(1+\d_N)^{l_{d+1}-l_1+\sum_{i=1}^{d-1}l_{i+1}}.
\]
Hence the total volume of this union is $\asymp_d \d_N^d(1+\d_N)^{l_{d+1}-l_1}$. Multiplication by the length $\asymp \d_N(1+\d_N)^{l_1}$ of the $n_1$-interval gives
\[
E(\mbf l,\signvec;N)\asymp_d \d_N^{d+1}(1+\d_N)^{l_{d+1}}.
\]
The equality of the two expected values in a paired sign choice is exact: the two target unions are sent to one another by a coordinate reflection on the torus, and the $n_1$-intervals have the same length after the first-coordinate sign is flipped.

To estimate the error for a fixed $\bma$, we rewrite the count on the interval for $n_1$ as a difference of two initial orbit sums. Beck's maximal discrepancy estimate \eqref{Thm of Beck}, applied with $\varphi(t)=t^2$, gives
\[
\#S_{\bma}(\mbf l,\signvec;N)
=E(\mbf l,\signvec;N)+O_{\bma}\big((\log (1+\d_N)^{l_1})^{d+1/2}\big).
\]
Since $l_1\ll \log N/\d_N$, the error is $O_{\bma}((\log N)^{d+1/2})$. On the other hand, the lower bound $(1+\d_N)^{l_{d+1}}\ge(\log N)^{3d+6}$ gives
\[
\d_N E(\mbf l,\signvec;N)\gg_d \d_N^{d+2}(\log N)^{3d+6}
=(\log N)^{d+2},
\]
which dominates $(\log N)^{d+1/2}$. This proves the claimed uniform estimate.

\end{proof}

\begin{remark}
The uniformity in the parameters appearing in the proposition above is worth making explicit. Beck's estimate \eqref{Thm of Beck} is a maximal discrepancy bound: for a fixed generic $\boldsymbol\alpha$, its constant is independent of the endpoint $\bf{x}$ of the anchored box and of the time parameter $N$. Each translated torus box occurring here, even if it crosses the boundary of the fundamental cube, is a signed sum of at most $O_d(1)$ anchored boxes. Hence the same full-measure set and the same implied constant, up to a factor depending only on $d$, apply simultaneously to all translations, all sign choices, all $\mathbf l\in L_2(N)$ and all $N$.
\end{remark}

For the sake of simplicity, we abbreviate $S_{\bma}(\mbf{l},\signvec;N)$ and denote it by $S(\mbf{l},\signvec;N)$ in later discussions. Using the lemma above, we can estimate the sum $\bar{D}_{5}$ by cancelling out the main terms.

\noindent\textbf{Proof of Proposition \ref{Proj two Prop cancellation of the main terms}}
\begin{proof}
Suppose that $\a$ satisfies Lemma \ref{Proj two d_N error estimation for the number of elements}. 
Denote a paired sign choice by $\mathfrak p=(\signvec^+,\signvec^-)$, where $\signvec^+$ and $\signvec^-$ are single sign vectors which differ in exactly one coordinate, and where the sign of the divisor $n_1\prod_{j=1}^{d} (n_1\a_j-n_{j+1})$ is positive on $S(\mbf l,\signvec^+;N)$ and negative on $S(\mbf l,\signvec^-;N)$. 
$$\bar{D}_{5}=\sum_{\substack{\mbf{l}\in L_2(N),\\ {\mathfrak p=(\signvec^+,\signvec^-)\in\mathcal P}}}
\l(\sum_{
\mbf{n}\in S (\mbf{l},\signvec^+;N)}
 \fc{g(\bm{n}, \bma; N)}{n_1\prod_{j=1}^{d} (n_1\a_j-n_{j+1})}+
 \sum_{
\mbf{n}\in S (\mbf{l},\signvec^-;N)}
\fc{g(\bm{n}, \bma; N)}{n_1\prod_{j=1}^{d} (n_1\a_j-n_{j+1})}\r).
$$
Inside each $\SlProjtwo$, write $g(\bm{n}, \bma; N)_{\max}$ and $g(\bm{n}, \bma; N)_{\min}$ for the maximal and minimal values of $g(\bm{n}, \bma; N)$ over $\mbf{n}\in \SlProjtwo$, then
$$
|g(\bm{n}, \bma; N)_{\max} - g(\bm{n}, \bma; N)_{\min}|\ll \d_N
$$
and 
$$
0\le g(\bm{n}, \bma; N)_{\min}\le g(\bm{n}, \bma; N)_{\max}\le 1.
$$
For a fixed paired sign choice $\mathfrak p$, the sum is estimated as follows, (the constants appearing in different lines of the inequalities may differ, but depend only on $\bma$ or $\varphi$.)
$$
\bal
&
\sum_{
\mbf{n}\in S (\mbf{l},\signvec^+;N)}
 \fc{g(\bm{n}, \bma; N)}{n_1\prod_{j=1}^{d} (n_1\a_j-n_{j+1})}+
 \sum_{
\mbf{n}\in S (\mbf{l},\signvec^-;N)}
\fc{g(\bm{n}, \bma; N)}{n_1\prod_{j=1}^{d} (n_1\a_j-n_{j+1})}\\
&
\le 
\fc{g(\bm{n}, \bma; N)_{\max}}{(1+\d_N)^{l_{d+1}}}(1+C\d_N)\ElProjtwo-\fc{g(\bm{n}, \bma; N)_{\min}}{(1+\d_N)^{l_{d+1}+d+1}}(1-C\d_N)\ElProjtwo\\
&\le
\fc{\l(g(\bm{n}, \bma; N)_{\max}(1+\d_N)^{d+1} -g(\bm{n}, \bma; N)_{\min}\r)}{(1+\d_N)^{l_{d+1}+d+1}}E(\mbf{l},\signvec;N)+C\fc{g(\bm{n}, \bma; N)_{\max}} {(1+\d_N)^{l_{d+1}}}\d_N\ElProjtwo\\
&\le
C\d_N \fc{\ElProjtwo}{(1+\d_N)^{l_{d+1}}} \\
&\le
C\d_N^{d+2}.\\
\eal
$$
The other direction is the same:
$$
\sum_{\mbf{n}\in S (\mbf{l},\signvec^+;N)}
 \fc{g(\bm{n}, \bma; N)}{n_1\prod_{j=1}^{d} (n_1\a_j-n_{j+1})}+
 \sum_{
\mbf{n}\in S (\mbf{l},\signvec^-;N)}
\fc{g(\bm{n}, \bma; N)}{n_1\prod_{j=1}^{d} (n_1\a_j-n_{j+1})}\\
\gg_{\bma}
-\d_N^{d+2}.
$$

Summing over all fixed pairs $\mathfrak p=(\signvec^+,\signvec^-)\in\mathcal P$ and over $l_1,\dots,l_{d+1}$, note that from \eqref{Proj two condition for l}, $0\le l_{j+1}\le l_1, \ 1\le j\le d$, we have:
$$
\bal
\l|\bar{D}_{5}\r|\ll_{\bma} \sum_{\mbf{l}\in L_2(N)} \d_N^{d+2}\ll_{\bma} \sum_{\fc{\log \log N }{\d_N}\le l_1 \le \fc{\log N}{\d_N}} \d_N^{d+2}l_1^{d}
\ll_{\bma} \d_N^{d+2}\cdot \fc{(\log N)^{d+1}}{\d_N^{d+1}}=(\log N)^{d-1}.
\eal
$$
\end{proof}

\section{Estimation of small exponentials}\label{section 5 estimation of small exponentials}
Finally we study the contribution of the linear forms $\L_\mbf{s}$ in \eqref{Proj two multiplied out term}.


 Throughout this section and in Appendix \ref{appendix:proof-of-proj-two-key-lemma}, $\mathfrak p=(\signvec^+,\signvec^-)\in\mathcal P$ always denotes one of the fixed paired sign choices introduced above; in particular the family $\mathcal P$ is a partition of the sign vectors and no sign cell is counted twice (see the definition before \eqref{Proj two condition for l}). 
Let 
\beq\label{Proj two exponential sums}
\bal
&\bar{D}_{6}=\bar{D}_{6}^{(\mbf{s})}=\sum_{\mbf{n}\in U_4(\bma;N)} \fc{e^{2\pi \i \L_{\mbf{s}}(\mbf{n})}}{n_1\prod_{j=1}^{d} (n_1\a_j-n_{j+1})}\cdot g(\bm{n}, \bma; N)\\
&=\sum_{\substack{\mbf{l}\in L_2(N),\\ {\mathfrak p=(\signvec^+,\signvec^-)\in\mathcal P}}}
\Bigg(\sum_{
\mbf{n}\in S (\mbf{l},\signvec^+;N)}
  \bar{f}(\mbf{n},x,\bma)+
 \sum_{
\mbf{n}\in S (\mbf{l},\signvec^-;N)}
  \bar{f}(\mbf{n},x,\bma),
\Bigg)\\
\eal
\eeq
where 
\beq\label{Proj two exponential term}
\bar{f}(\mbf{n},x,\bma)=\fc{e^{2\pi \i \L_{\mbf{s}}(\mbf{n})}}{n_1\prod_{j=1}^{d} (n_1\a_j-n_{j+1})}\cdot g(\bm{n}, \bma; N)
\eeq
and $\L$ is one of the $2^{d+1}-1$ linear forms $\L_{\mbf{s},x,N, \bma}$ defined in \eqref{Proj two linear forms}.

We now define bad vectors which encode the cells whose exponential sums do not cancel out. 
\begin{definition}\label{Proj two big vector definition}
For a paired sign choice $\mathfrak p=(\signvec^+,\signvec^-)$, we say that $\mbf l=(l_1,\dots,l_{d+1})$ is \textit{$\mathfrak p$-big} if
\beq\label{Proj two big vector}
\frac{|S(\mbf l,\signvec^+;N)|+|S(\mbf l,\signvec^-;N)|}{\log N}\le
\left|\sum_{\mbf n\in S(\mbf l,\signvec^+;N)} e^{2\pi \i \L(\mbf n)}-\sum_{\mbf n\in S(\mbf l,\signvec^-;N)} e^{2\pi \i \L(\mbf n)}\right|.
\eeq
\end{definition}

We need to control the number of bad vectors defined above. The following lemma is an adaptation of Beck's ``Key Lemma'' to our setting; the proof is deferred to Appendix \ref{appendix:proof-of-proj-two-key-lemma}.

\begin{lemma}\label{Proj two number of e-big vectors}
For almost every $\bma$, uniformly for every $N$, every $x\in[0,1]$, every nonzero $\mbf{s}\in\{0,1\}^{d+1}$, and every paired sign choice $\mathfrak p$, the number of $\mathfrak p$-big vectors associated with $\L_{\mbf{s},x,N,\bma}$ is
$$
\ll_{\bma}(\log \log N)(\log N)^d\d_N^{-(d+1)}.
$$
\end{lemma}
\begin{remark}
A preview of the proof is that, instead of considering one vector at a time, for each fixed exponential phase $\L_{\mbf{s}}$, we will consider at once \emph{a set of vectors} $\mbf l$ in a special arithemetic progression. These sets are called \emph{special lines} (see Definition \ref{Proj two defnition for e-line}). Lemma \ref{Proj two key lemma} shows that each such line contains \emph{at most one} $\mathfrak {p}$-big vector, whose corresponding paired cells without sufficient cancellation, as in \eqref{Proj two big vector}. Finally, the number of possible special lines is much smaller since they contains arithemtic progressions of vectors. This gives the desired logarithmic saving.
\end{remark}

The next lemma states that, inside a signed cell, the weight $w(\mbf n)$ is nearly constant; on the two cells of a paired sign choice the leading constants have opposite signs. Hence the cancellation in the unweighted exponential sums transfers to the weighted sums.
\begin{lemma}\label{Proj two weighted cancellation lemma}
Fix a paired sign choice $\mathfrak p=(\signvec^+,\signvec^-)$ as in Definition \ref{Proj two big vector definition}, and assume that $\mbf l$ is not $\mathfrak p$-big. Put
\[
w(\mbf n)=\frac{g(\mbf n,\bma;N)}{n_1\prod_{j=1}^{d}(n_1\alpha_j-n_{j+1})}.
\]
Then
\[
\begin{aligned}
&\left|\sum_{\mbf n\in S(\mbf l,\signvec^+;N)} e^{2\pi \i \L(\mbf n)}w(\mbf n)
+\sum_{\mbf n\in S(\mbf l,\signvec^-;N)} e^{2\pi \i \L(\mbf n)}w(\mbf n)\right|  \\
&\qquad\qquad\ll
(1+\delta_N)^{-l_{d+1}}\frac{E(\mbf l,\signvec^+;N)+E(\mbf l,\signvec^-;N)}{\log N}.
\end{aligned}
\]
Equivalently, by Lemma \ref{Proj two d_N error estimation for the number of elements}, the right-hand side is
$\ll (1+\delta_N)^{-l_{d+1}}E(\mbf l,\signvec^+;N)/\log N$.
\end{lemma}

\begin{proof}
For $\mbf n\in S(\mbf l,\signvec^{\pm};N)$ the absolute divisor satisfies
\[
|n_1|\prod_{j=1}^{d}|n_1\alpha_j-n_{j+1}|=(1+O(\delta_N))(1+\delta_N)^{l_{d+1}}.
\]
Moreover, because $|n_1|/N^d\le1/4$ on $U_4$ and $|n_1\alpha_j-n_{j+1}|<1/2$, the function
$t\mapsto(\sin(2\pi t)/(2\pi t))^2$ is Lipschitz on all intervals involved. Hence there exists a number $\omega_{\mbf l}$ with
$|\omega_{\mbf l}|\ll(1+\delta_N)^{-l_{d+1}}$ such that
\[
w(\mbf n)=\omega_{\mbf l}+O\bigl(\delta_N(1+\delta_N)^{-l_{d+1}}\bigr),\qquad
\mbf n\in S(\mbf l,\signvec^+;N),
\]
and
\[
w(\mbf n)=-\omega_{\mbf l}+O\bigl(\delta_N(1+\delta_N)^{-l_{d+1}}\bigr),\qquad
\mbf n\in S(\mbf l,\signvec^-;N).
\]
Therefore
\[
\begin{aligned}
&\left|\sum_{S(\mbf l,\signvec^+;N)}e^{2\pi \i \L(\mbf n)}w(\mbf n)+\sum_{S(\mbf l,\signvec^-;N)}e^{2\pi \i \L(\mbf n)}w(\mbf n)\right| \\
&\quad\ll (1+\delta_N)^{-l_{d+1}}
\left|\sum_{S(\mbf l,\signvec^+;N)}e^{2\pi \i \L(\mbf n)}-\sum_{S(\mbf l,\signvec^-;N)}e^{2\pi \i \L(\mbf n)}\right| \\
&\qquad+\delta_N(1+\delta_N)^{-l_{d+1}}\bigl(|S(\mbf l,\signvec^+;N)|+|S(\mbf l,\signvec^-;N)|\bigr).
\end{aligned}
\]
Since $\mbf l$ is not $\mathfrak p$-big, the first term is bounded by
$(1+\delta_N)^{-l_{d+1}}(|S(\mbf l,\signvec^+;N)|+|S(\mbf l,\signvec^-;N)|)/\log N$. Lemma \ref{Proj two d_N error estimation for the number of elements} and the equality of the paired expected values give
\[
|S(\mbf l,\signvec^+;N)|+|S(\mbf l,\signvec^-;N)|\ll_{\bma}E(\mbf l,\signvec^+;N)+E(\mbf l,\signvec^-;N).
\]
Since $\delta_N=(\log N)^{-2}$, the error term is smaller than the first bound. The claim follows.
\end{proof}
Using Lemma \ref{Proj two number of e-big vectors}, we are ready to estimate the contribution of the exponential terms. We prove the following proposition:
\begin{proposition}\label{Proj two Estimation of the exponential terms}
For almost every $\bma$, we have 
$$
\l|\bar{D}_{6}\r|\ll (\log N)^d(\log \log N)
$$
\end{proposition}
\begin{proof}
Suppose $\boldsymbol{\alpha}$ satisfies Lemma \ref{Proj two number of e-big vectors}. We consider a fixed
$\bar D_6^{(\boldsymbol{s})}$; the number of possible nonzero
$\boldsymbol{s}\in\{0,1\}^{d+1}$ is finite and depends only on $d$.  Decompose
\[
\bar D_6^{(\boldsymbol{s})}=\sum_{\mathrm{small}}+\sum_{\mathrm{big}},
\]
where
\[
\sum_{\mathrm{small}}
=\sum_{\mathfrak p=(\signvec^+,\signvec^-)\in\mathcal P}
\sum_{\substack{\boldsymbol l\in L_2(N)\\ \boldsymbol l\text{ not }\mathfrak p\text{-big}}}
\left(
\sum_{\boldsymbol n\in S(\boldsymbol l,\signvec^+;N)}\bar f(\boldsymbol n,x,\boldsymbol\alpha)+\sum_{\boldsymbol n\in S(\boldsymbol l,\signvec^-;N)}\bar f(\boldsymbol n,x,\boldsymbol\alpha)
\right),
\]
and $\sum_{\mathrm{big}}$ is defined similarly with $\boldsymbol l$ being $\mathfrak p$-big.
For the small vectors, Lemma \ref{Proj two weighted cancellation lemma}, Lemma
\ref{Proj two d_N error estimation for the number of elements}, and the range
\eqref{Proj two condition for l} give
\[
\begin{aligned}
\left|\sum_{\mathrm{small}}\right|
&\ll
\sum_{\boldsymbol l\in L_2(N)}(1+\delta_N)^{-l_{d+1}}
\frac{E(\boldsymbol l,\signvec^+;N)+E(\boldsymbol l,\signvec^-;N)}{\log N} \\
&\ll
\sum_{\boldsymbol l\in L_2(N)}\frac{\delta_N^{d+1}}{\log N}
\ll
\frac{\delta_N^{d+1}}{\log N}\frac{(\log N)^{d+1}}{\delta_N^{d+1}}
\ll(\log N)^d.
\end{aligned}
\]
For the big vectors, Lemma \ref{Proj two number of e-big vectors} and
$E(\boldsymbol l,\signvec^{\pm};N)\asymp_d\delta_N^{d+1}(1+\delta_N)^{l_{d+1}}$ imply
\[
\begin{aligned}
\left|\sum_{\mathrm{big}}\right|
&\ll
\sum_{\boldsymbol l\text{ is }\mathfrak p\text{-big}}
(1+\delta_N)^{-l_{d+1}}\bigl(E(\boldsymbol l,\signvec^+;N)+E(\boldsymbol l,\signvec^-;N)\bigr)\\
&\ll
\bigl((\log\log N)(\log N)^d\delta_N^{-(d+1)}\bigr)\delta_N^{d+1}
\ll(\log N)^d\log\log N.
\end{aligned}
\]
Summing over the finitely many choices of $\boldsymbol{s}$ proves the proposition.
\end{proof}
Combining Proposition \ref{Proj two Prop Dbar}, \ref{Proj two D4}, \ref{Proj two Prop cancellation of the main terms} and \ref{Proj two Estimation of the exponential terms}, we complete the proof of the convergent part \eqref{Proj two convergent part of the theorem} of Theorem \ref{Proj two main result}.
\vspace{+0.3cm}
\section{Divergent part of the main theorem}\label{divergent part of the main theorem}
\noindent \begin{proof}[\textbf{Proof of the divergent part \eqref{Proj two divergent part of the theorem} of Theorem \ref{Proj two main result}}]
We prove the divergent part using the anchored discrepancy $[0,x)$ and the Fourier expansion \eqref{Final expression of the smoothed discrepancy}. This avoids the sign-reflected frequency which appears if one projects the symmetric interval $[-x,x]$ onto sine functions as in \cite{Beck}.

By Gallagher's metric theorem for simultaneous Diophantine approximation \cite{gallagher1962metric}, we have the following lemma.

\begin{lemma}\label{Gallagher lower bound in divergent part}
Let $\varphi(n)$ be positive and nondecreasing, and suppose that $\sum_{n=1}^{\infty}1/\varphi(n)=\infty$. Then for almost every $\bma=(\alpha_1,\dots,\alpha_d)\in\mathbb{R}^d$, there exist infinitely many $n\in\mathbb{N}$ such that
\beq\label{Proj two Gallagher small product}
n\prod_{j=1}^{d}\|n\alpha_j\|<\frac{1}{(\log n)^d\varphi(\log\log n)}.
\eeq
\end{lemma}

Fix $0<\varepsilon<1/2$. Since $\varphi$ is nondecreasing and $\sum_n1/\varphi(n)=\infty$, the reciprocal sum for $t\mapsto\varphi(2t)$ also diverges. We apply Lemma \ref{Gallagher lower bound in divergent part} with $t\mapsto\varphi(2t)$ in place of $\varphi(t)$. 
Take $\bma$ in the full-measure set satisfying this Gallagher conclusion and the corresponding lower-bound estimates in dimensions $1,d-1$, and $d$.  More precisely, we use the same Borel--Cantelli argument as in Lemma \ref{Proj two absolute lower bound the the divisor}, with thresholds $1/(n(\log n)^{r+\varepsilon})$ in dimensions $r=1,d-1,d$.
Let $n_1^*$ be one of the infinitely many integers thereby obtained, and let $n_{j+1}^*$ be the nearest integer to $n_1^*\alpha_j$ for $1\le j\le d$. Then, for $n_1^*$ sufficiently large,
\beq\label{Proj two product lower and upper in divergent part}
\frac{1}{n_1^*(\log n_1^*)^{d+\varepsilon}}
\le \prod_{j=1}^{d}\|n_1^*\alpha_j-n_{j+1}^*\|
\le \frac{1}{n_1^*(\log n_1^*)^d\varphi(2\log\log n_1^*)}.
\eeq
Moreover, for every $1\le m\le d$,
\beq\label{Proj two partial product lower divergent part}
\prod_{j\neq m}\|n_1^*\alpha_j-n_{j+1}^*\|\gg_{\bma,\varepsilon}\frac{1}{n_1^*(\log n_1^*)^{d-1+\varepsilon}},
\eeq
and combining this with the upper bound in \eqref{Proj two product lower and upper in divergent part} gives
\beq\label{Proj two individual lower upper divergent part}
\frac{1}{n_1^*(\log n_1^*)^{1+\varepsilon}}\ll_{\bma,\varepsilon}\|n_1^*\alpha_j-n_{j+1}^*\|\ll \frac{1}{(\log n_1^*)^{1/2}},
\qquad 1\le j\le d.
\eeq

We next choose $N^*$ so that the oscillating numerators in the resonant Fourier mode are not small.

\begin{lemma}\label{Proj two choose N in divergent part}
For all sufficiently large $n_1^*$, there exists an integer
$$
N^*\in\left[
\frac{1}{\prod_{j=1}^{d}\|n_1^*\alpha_j-n_{j+1}^*\|},
\frac{2}{\prod_{j=1}^{d}\|n_1^*\alpha_j-n_{j+1}^*\|}\right]
$$
such that
$$
\|N^*(n_1^*\alpha_j-n_{j+1}^*)\|\ge \frac{1}{16d},\qquad 1\le j\le d.
$$
Consequently,
\beq\label{Proj two numerator lower divergent part}
\prod_{j=1}^{d}\bigl|1-e^{2\pi iN^*(n_1^*\alpha_j-n_{j+1}^*)}\bigr|\gg_d 1.
\eeq
\end{lemma}

\begin{proof}
Let
$$
I=\left[
\frac{1}{\prod_{j=1}^{d}\|n_1^*\alpha_j-n_{j+1}^*\|},
\frac{2}{\prod_{j=1}^{d}\|n_1^*\alpha_j-n_{j+1}^*\|}\right].
$$
For fixed $j$, define
$$
B_j=\{m\in I\cap\mathbb{Z}: \|m(n_1^*\alpha_j-n_{j+1}^*)\|<1/(16d)\}.
$$
We use the elementary fact that, for any interval $J$ and any $0<\theta<1/2$,
\[
\#\{m\in J\cap\mathbb Z:\|m\theta\|<\rho\}
\le 2\rho\#(J\cap\mathbb Z)+O(\theta^{-1}+1),
\]
which follows by splitting $J$ into blocks of length $\lfloor\theta^{-1}\rfloor$ and using periodicity modulo one.
In an interval of length comparable to $1/\prod_{j=1}^{d}\|n_1^*\alpha_j-n_{j+1}^*\|$, the number of such integers is bounded by
$$
\#B_j\le \frac{1}{8d}\#(I\cap\mathbb{Z})+O\left(\frac{1}{\|n_1^*\alpha_j-n_{j+1}^*\|}+1\right).
$$
Since
$$
\frac{1/\|n_1^*\alpha_j-n_{j+1}^*\|}{1/\prod_{i=1}^{d}\|n_1^*\alpha_i-n_{i+1}^*\|}
=\prod_{i\neq j}\|n_1^*\alpha_i-n_{i+1}^*\|\le(\log n_1^*)^{-(d-1)/2}=o(1),
$$
we have $\#B_j\le(4d)^{-1}\#(I\cap\mathbb{Z})$ for large $n_1^*$. Therefore $\#\cup_{j=1}^{d}B_j<\#(I\cap\mathbb{Z})$, and an admissible $N^*$ exists. The last assertion follows from $|1-e^{2\pi it}|\gg\|t\|$.
\end{proof}

We now project the anchored smoothed discrepancy onto the complex frequency $n_1^*$:
\beq\label{Proj two anchored projection divergent part}
I_{n_1^*}:=\int_0^1\bar D(\bma,x;N^*)e^{-2\pi i n_1^*x}\,dx.
\eeq
In the Fourier expansion \eqref{Final expression of the smoothed discrepancy}, the contribution with $n_1=0$ is zero: if some $n_{j+1}\neq0$, then the numerator $1-e^{2\pi i N^*n_{j+1}}$ vanishes, while the all-zero mode has already been cancelled by the volume term. For $n_1\neq0$, orthogonality gives
$$
\int_0^1(1-e^{2\pi i n_1x})e^{-2\pi i n_1^*x}\,dx=-\mathbf 1_{\{n_1=n_1^*\}}.
$$
Then
\beq\label{Proj two I equals S1 S2}
I_{n_1^*}=S_1+S_2,
\eeq
where the main term is
\beq\label{Proj two S1 displayed expression}
\begin{aligned}
S_1
&=-\frac{i^{d+1}}{2\pi n_1^*}\left(\frac{\sin(2\pi n_1^*/(N^*)^d)}{2\pi n_1^*/(N^*)^d}\right)^2\\
&\times\prod_{j=1}^{d}\left(
\frac{1-e^{-2\pi iN^*(n_1^*\alpha_j-n_{j+1}^*)}}{2\pi(n_1^*\alpha_j-n_{j+1}^*)}
\left(\frac{\sin(2\pi(n_1^*\alpha_j-n_{j+1}^*))}{2\pi(n_1^*\alpha_j-n_{j+1}^*)}\right)^2
\right),
\end{aligned}
\eeq
and the remainder with the same first coordinate is
\beq\label{Proj two S2 displayed expression}
\begin{aligned}
S_2
&=-\frac{i^{d+1}}{2\pi n_1^*}\left(\frac{\sin(2\pi n_1^*/(N^*)^d)}{2\pi n_1^*/(N^*)^d}\right)^2\\
&\times\sum_{\mathbf m=(m_1,\dots,m_d)\in\mathbb{Z}^d\setminus\{\mathbf m^*\}}\prod_{j=1}^{d}\left(
\frac{1-e^{-2\pi iN^*(n_1^*\alpha_j-m_j)}}{2\pi(n_1^*\alpha_j-m_j)}
\left(\frac{\sin(2\pi(n_1^*\alpha_j-m_j))}{2\pi(n_1^*\alpha_j-m_j)}\right)^2
\right).
\end{aligned}
\eeq

From \eqref{Proj two product lower and upper in divergent part} and \eqref{Proj two numerator lower divergent part}, and since all smoothing factors in \eqref{Proj two S1 displayed expression} are bounded below by an absolute positive constant for large $n_1^*$, we obtain
\beq\label{Proj two S1 lower divergent part}
|S_1|\gg_{\bma,d}\frac{1}{n_1^*\prod_{j=1}^{d}\|n_1^*\alpha_j-n_{j+1}^*\|}\gg_{\bma,d}(\log n_1^*)^d\varphi(2\log\log n_1^*).
\eeq
The lower and upper bounds in \eqref{Proj two product lower and upper in divergent part} imply, respectively,
\[
N^*\ll n_1^*(\log n_1^*)^{d+\varepsilon},
\qquad
N^*\ge n_1^*(\log n_1^*)^d\varphi(2\log\log n_1^*).
\]
In particular $N^*\to\infty$ along the chosen sequence and
\[
\log N^*\asymp \log n_1^*,
\qquad
\log\log N^*\le 2\log\log n_1^*.
\]
for all sufficiently large $n_1^*$. Therefore, by monotonicity of $\varphi$,
\beq\label{Proj two S1 in terms of Nstar}
|S_1|\gg_{\bma,\varphi}(\log N^*)^d\varphi(\log\log N^*).
\eeq

It remains to estimate $S_2$. For a non-main term put
\[
J=\{j:m_j=n_{j+1}^*\}.
\]
Since the term is not the main term, $J\ne\{1,\ldots,d\}$. If $j\in J$, then the corresponding denominator contributes the small factor $\|n_1^*\alpha_j-n_{j+1}^*\|^{-1}$. If $j\notin J$, then $m_j$ is not the nearest integer $n_{j+1}^*$; hence $|n_1^*\alpha_j-m_j|\ge1/2$. As in \eqref{inequality to ignore n such that two or more factors are big}, the factor
\[
\frac{1}{|n_1^*\alpha_j-m_j|}
\left|\frac{\sin 2\pi(n_1^*\alpha_j-m_j)}{2\pi(n_1^*\alpha_j-m_j)}\right|^2
\]
is $O(r^{-3})$ for $r\le |n_1^*\alpha_j-m_j|<r+1$ and $r=1,2,\ldots.$
Thus every free non-nearest coordinate contributes a bounded factor $\sum_{r\ge1}r^{-3}=O(1)$. Consequently
\beq\label{Proj two S2 upper divergent part}
\begin{aligned}
|S_2|
&\ll \sum_{J\subsetneq\{1,\dots,d\}}
\frac{1}{n_1^*}\prod_{j\in J}\frac{1}{\|n_1^*\alpha_j-n_{j+1}^*\|} \\
&\ll \sum_{m=1}^{d}
\frac{1}{n_1^*\prod_{j\neq m}\|n_1^*\alpha_j-n_{j+1}^*\|}+O(1) \\
&\ll_{\bma,\varepsilon}(\log n_1^*)^{d-1+\varepsilon},
\end{aligned}
\eeq
where the last step uses the $(d-1)$-fold product lower bounds in \eqref{Proj two partial product lower divergent part}.

This is negligible compared with \eqref{Proj two S1 in terms of Nstar}. Hence, for infinitely many $N^*$,
\[
|I_{n_1^*}|\le\int_0^1|\bar D(\bma,x;N^*)|\,dx
\le\sup_{0<x\le1}|\bar D(\bma,x;N^*)|,
\]
and therefore
\[
\sup_{0<x\le1}|\bar D(\bma,x;N^*)|
\ge |I_{n_1^*}|\gg_{\bma,\varphi}(\log N^*)^d\varphi(\log\log N^*).
\]
Finally, Proposition \ref{Proj two Prop Dbar} gives the same lower bound for the original anchored discrepancy, since the smoothing error is $O_{\bma,\varepsilon}((\log N^*)^{d-1+\varepsilon})$, which is negligible. Therefore
$$
\Delta(\bma;N^*)\gg_{\bma,\varphi}(\log N^*)^d\varphi(\log\log N^*)
$$
for infinitely many $N^*$, proving the divergent claim.
\end{proof}
\section*{Acknowledgements}

This research was funded in part by the Austrian Science Fund (FWF) [10.55776/PAT5120424]. For open access purposes, the author has applied a CC BY public copyright license to any author accepted manuscript version arising from this submission.

\appendix
\section{Proof of Lemma \ref{Proj two lem1}}\label{appendix:proof-of-proj-two-lem1}

\begin{proof}[Proof of Lemma \ref{Proj two lem1}]
Note that the integral 
$$
J(n)=\int_{[0,1]^d} \fc{1}{\prod_{j=1}^{d} \l(\|n \a_j\|\varphi(|\log \|n\a_j\||)\r)} d\bma
$$
is finite and independent of $n$. 
Indeed, near zero
\[
\int_0^{1/2}\frac{dt}{t\varphi(|\log t|)}=
\int_{\log 2}^{\infty}\frac{du}{\varphi(u)}<\infty,
\]
where the last inequality follows from the monotonicity of $\varphi$ and the convergence of $\sum 1/\varphi(n)$. Also the series $\sum\limits_{n=2}^\infty \fc{1}{n\varphi(\log n)}$ converges, so the integral
\beq{\label{Proj two int-series1}}
\int_{[0,1]^d} \l(\sum_{n=2}^\infty \fc{1}{n\varphi(\log n)\prod_{j=1}^{d} \l(\|n \a_j\|\varphi(|\log \|n\a_j\||)\r)}\r)d\bma=\sum_{n=2}^\infty \fc{J(n)}{n\varphi(\log n)}
\eeq
is finite.
Therefore we have that the series 
 $$\sum_{n=2}^\infty \fc{1}{n\varphi(\log n)\prod_{j=1}^{d} \l(\|n \a_j\|\varphi(|\log \|n\a_j\||)\r)}$$
converges for almost every $\bma\in [0,1]^d$. The periodicity of $\|\cdot\|$ gives the result for almost every $\bma\in \R^d $.
\end{proof}

\section{Proof of Lemma \ref{Proj two absolute lower bound the the divisor}}\label{appendix:proof-of-proj-two-absolute-lower-bound}
\begin{proof}[Proof of Lemma \ref{Proj two absolute lower bound the the divisor}]
Let $\psi$ be as in Lemma \ref{Proj two absolute lower bound the the divisor}.  For $|n|$ large, define
\[
E_n:=\left\{\bma\in[0,1]^d:
|n|\prod_{j=1}^{d}\|n\alpha_j\|<
\frac{1}{(\log |n|)^d\psi(\log\log |n|)}
\right\}.
\]
The standard multiplicative-cusp estimate gives
\[
\operatorname{Leb}(E_n)
\ll_d
\frac{1}{|n|(\log |n|)\psi(\log\log |n|)}.
\]
Hence $\sum_{n\in\mathbb Z\setminus\{0\}}\operatorname{Leb}(E_n)<\infty$, and the Borel--Cantelli lemma implies that, for almost every $\bma$,
\[
|n|\prod_{j=1}^{d}\|n\alpha_j\|
\gg_{\bma,\psi}
\frac{1}{(\log |n|)^d\psi(\log\log |n|)}.
\]
The one-dimensional estimate for each $\|n\alpha_j\|$ follows from the same argument applied to the sets
\[
\left\{\alpha_j\in[0,1]:\|n\alpha_j\|<\frac{1}{|n|\psi(\log |n|)}\right\}.
\]
\end{proof}

\section{A self-contained weaker substitute for Lemma \ref{fregoli-reciprocal-input}}\label{appendix:weak-fregoli-substitute}

The proof of Proposition \ref{control for sum when one coordinate is large} uses the sharp estimate in \cite{fregoli2024sumsreciprocalsfractionalparts}. The following elementary substitute is included only to provide a short self-contained alternative at the cost of extra powers of $\log\log N$.
\begin{proposition}\label{weak-fregoli-substitute}
Let $r\ge1$ and let $\varrho$ be a positive nondecreasing function with $\sum_m1/\varrho(m)<\infty$. For almost every $\bm\beta\in\mathbb R^r$,
\[
        \sum_{1\le n\le X}\frac{1}{n\prod_{i=1}^r\|n\beta_i\|}
        \ll_{\bm\beta,r,\varrho}(\log X)^{r+1}\varrho^{r+1}(\log\log X),\qquad X\ge30.
\]
In particular, when $r=d-1$ and $X=N^d(\log N)^d$, the bound needed in Proposition \ref{control for sum when one coordinate is large} becomes
\[
        (\log N)^d\varrho^d(\log\log N).
\]
\end{proposition}
\begin{proof}
Apply Lemma \ref{Proj two lem1} with the auxiliary function
\[
        \varphi_0(t)=t\varrho(\log t).
\]
Then for almost every $\bm\beta$ the series
\[
\sum_{n=3}^{\infty}
\frac{1}{n(\log n)\varrho(\log\log n)
\prod_{i=1}^{r}\bigl(\|n\beta_i\|\,|\log\|n\beta_i\||\,\varrho(\log |\log\|n\beta_i\||)\bigr)}
\]
converges. The elementary Borel--Cantelli estimate $\|n\beta_i\|\ge n^{-2}$ for all sufficiently large $n$ gives $|\log\|n\beta_i\||\le2\log n$. Hence, for $n\le X$,
\[
\begin{aligned}
\frac{1}{n\prod_i\|n\beta_i\|}
&\ll_{\bm\beta,r,\varrho}(\log X)^{r+1}\varrho^{r+1}(\log\log X) \\
&\quad\times
\frac{1}{n(\log n)\varrho(\log\log n)
\prod_i\bigl(\|n\beta_i\|\,|\log\|n\beta_i\||\,\varrho(\log |\log\|n\beta_i\||)\bigr)}.
\end{aligned}
\]
Summing over $n\le X$ and using the convergence of the last series proves the claim.
\end{proof}

\section{Proof of Lemma \ref{Proj two number of e-big vectors}}\label{appendix:proof-of-proj-two-key-lemma}

Throughout this appendix we fix $N$, a paired sign choice
$\mathfrak p=(\signvec^+,\signvec^-)$, and a linear form
$\L=\L_{\mbf{s},x,N,\bma}$. After Lemma \ref{Proj two d_N error estimation for the number of elements} has been established, no further measure-theoretic exceptional set is introduced in this appendix. The argument below is uniform for every phase $\L$; consequently the same full-measure set of parameters works uniformly for every $x$, every nonzero vector $\mbf s$, and every paired sign choice $\mathfrak p$.

Put
\[
A_N=\left\lfloor\frac{12\log\log N}{\log(1+\delta_N)}\right\rfloor,\qquad
R_N=(1+\delta_N)^{A_N}.
\]
Then $R_N\asymp(\log N)^{12}$. 
Recall that for a single sign vector $\signvec\in\{\pm1\}^{d+1}$, the set $S(\mbf l,\signvec;N)$ is described by the signed coordinates:
\[
\epsilon_1n_1\in I_1(\mbf l),\qquad
\epsilon_{j+1}(n_1\alpha_j-n_{j+1})\in I_{j+1}(\mbf l),\quad 1\le j\le d,
\]
where
\[
I_1(\mbf l)=[(1+\delta_N)^{l_1},(1+\delta_N)^{l_1+1}),
\]
\[
I_{j+1}(\mbf l)=[(1+\delta_N)^{-l_{j+1}},(1+\delta_N)^{-l_{j+1}+1}),\qquad 1\le j\le d-1,
\]
and
\[
I_{d+1}(\mbf l)=[(1+\delta_N)^{\kappa(\mbf l)},(1+\delta_N)^{\kappa(\mbf l)+1}),
\]
together with
\[
\kappa(\mbf l)=l_{d+1}-l_1+\sum_{i=1}^{d-1}l_{i+1}.
\]
Note that the closest-integer condition in the definition of $S(\mbf l,\signvec;N)$ is preserved by the imposed range.

\begin{definition}\label{Proj two defnition for e-line}
Let
\[
\mbf s_N=(A_N,-A_N,\dots,-A_N,(d+1)A_N)\in\mathbb Z^{d+1},
\]
where the entry $-A_N$ occurs in the positions $2,\dots,d$. Two vectors $\mbf l,\mbf h\in L_2(N)$ are called neighbors, written $\mbf l\to\mbf h$, if
\[
\mbf h=\mbf l+\mbf s_N.
\]
A special line is a maximal set of the form
\[
\{\mbf l_0+r\mbf s_N:r\in\mathbb Z\}\cap L_2(N).
\]
\end{definition}
If $\mbf h=\mbf l+\mbf s_N$, then
\beq\label{neighbor condition with RN}
I_i(\mbf h)=R_N I_i(\mbf l)
\quad(1\le i\le d+1)
\eeq
in the multiplicative sense: every endpoint of $I_i(\mbf l)$ is multiplied by $R_N$.

\begin{lemma}\label{Proj two key lemma}
On the full-measure set of $\bma$ given by Lemma \ref{Proj two d_N error estimation for the number of elements}, uniformly in $x$ and $\mbf s$, every special line contains at most one $\mathfrak p$-big vector.
\end{lemma}

\begin{proof}
Assume that a special line contains two $\mathfrak p$-big vectors $\mbf h^{(p)}$ and $\mbf h^{(q)}$, with $p<q$. We first prove that there exists
\[
\mbf n^*\in S(\mbf h^{(p)},\signvec^+;N)\cup S(\mbf h^{(p)},\signvec^-;N)
\]
such that
\beq\label{Proj two big L(n)}
\|\L(\mbf n^*)\|>(\log N)^{-2}.
\eeq
If not, then $e^{2\pi \i\L(\mbf n)}=1+O((\log N)^{-2})$ on both paired cells. By Lemma \ref{Proj two d_N error estimation for the number of elements} and by
$E(\mbf h^{(p)},\signvec^+;N)=E(\mbf h^{(p)},\signvec^-;N)$,
\[
\begin{aligned}
&\left|\sum_{S(\mbf h^{(p)},\signvec^+;N)}e^{2\pi \i\L(\mbf n)}-\sum_{S(\mbf h^{(p)},\signvec^-;N)}e^{2\pi \i\L(\mbf n)}\right|\\
&\qquad\ll_{\bma}
\frac{|S(\mbf h^{(p)},\signvec^+;N)|+|S(\mbf h^{(p)},\signvec^-;N)|}{(\log N)^2},
\end{aligned}
\]
contradicting the $\mathfrak p$-big condition. This proves \eqref{Proj two big L(n)}.

For a single sign vector $\signvec\in\{\signvec^+,\signvec^-\}$, define the inner part
$S_{\rm in}(\mbf h^{(q)},\signvec;N)$ by replacing each interval
$[A,(1+\delta_N)A)$ in the definition of $S(\mbf h^{(q)},\signvec;N)$ by
\beq\label{definition for inner points}
[(1+\delta^2_N)A,(1-\delta^2_N)(1+\delta_N)A).
\eeq
Let
\[
S_{\rm bd}(\mbf h^{(q)},\signvec;N)=S(\mbf h^{(q)},\signvec;N)\setminus S_{\rm in}(\mbf h^{(q)},\signvec;N).
\]

We now check stability under the progression direction $\mbf n^*$. If
$\mbf m\in S_{\rm in}(\mbf h^{(q)},\signvec;N)$ and $|r|\le(\log N)^4$, then the first signed coordinate changes by at most
\[
|r n_1^*|\le(\log N)^4(1+\delta_N)^{h^{(p)}_1+1}\ll(\log N)^4R_N^{-1}(1+\delta_N)^{h^{(q)}_1}=o(\delta^2_N(1+\delta_N)^{h^{(q)}_1}),
\]

Here we use the multiplicative endpoint relation
\[
(1+\delta_N)^{h^{(p)}_1+1}\le (1+\delta_N)R_N^{-1}(1+\delta_N)^{h^{(q)}_1},
\]
which follows from \eqref{neighbor condition with RN}.

Recall that $R_N\asymp(\log N)^{12}$ and $\delta_N=(\log N)^{-2}$, thus for $1\le j\le d-1$,
\[
|r(n_1^*\alpha_j-n_{j+1}^*)|\le(\log N)^4(1+\delta_N)^{-h^{(p)}_{j+1}+1}\ll(\log N)^4R_N^{-1}(1+\delta_N)^{-h^{(q)}_{j+1}}=o(\delta^2_N(1+\delta_N)^{-h^{(q)}_{j+1}}),
\]
and for the last coordinate,
\[
|r(n_1^*\alpha_d-n_{d+1}^*)|\le(\log N)^4(1+\delta_N)^{\kappa(\mbf h^{(p)})+1}\ll(\log N)^4R_N^{-1}(1+\delta_N)^{\kappa(\mbf h^{(q)})}=o(\delta^2_N(1+\delta_N)^{\kappa(\mbf h^{(q)})}).
\]
Hence
\beq\label{Proj two progression stability refined}
\mbf m+r\mbf n^*\in S(\mbf h^{(q)},\signvec;N),\qquad |r|\le(\log N)^4.
\eeq

For each coset of $\mathbb Z\mbf n^*$ in $\mathbb Z^{d+1}$, intersect the coset with $S(\mbf h^{(q)},\signvec;N)$ and decompose the result into maximal consecutive arithmetic progressions in the direction $\mbf n^*$. These progressions are disjoint and cover the whole cell. If a maximal progression meets the inner part, then \eqref{Proj two progression stability refined} shows that it contains at least $2(\log N)^4+1$ consecutive points, unless the missing consecutive points are outside the cell only because the progression has reached the boundary strip. All such endpoint pieces are therefore contained in $S_{\rm bd}(\mbf h^{(q)},\signvec;N)$ and will be estimated separately. Hence every remaining long progression has length $L_P\ge(\log N)^4$. If $P=\{\mbf m_0+r\mbf n^*:0\le r<L_P\}$ is one of these long progressions, then using \eqref{Proj two big L(n)},

\[
\left|\sum_{\mbf n\in P}e^{2\pi \i\L(\mbf n)}\right|
=\left|e^{2\pi \i\L(\mbf m_0)}\sum_{r=0}^{L_P-1}e^{2\pi \i r\L(\mbf n^*)}\right|
\ll\frac1{\|\L(\mbf n^*)\|}\ll(\log N)^2\ll\frac{L_P}{(\log N)^2}.
\]
Summing over all progressions gives
\beq\label{Proj two AP estimate revised}
\left|\sum_{\mbf n\in S(\mbf h^{(q)},\signvec;N)}e^{2\pi \i\L(\mbf n)}\right|
\ll
\frac{|S(\mbf h^{(q)},\signvec;N)|}{(\log N)^2}+|S_{\rm bd}(\mbf h^{(q)},\signvec;N)|.
\eeq

It remains to estimate the boundary part $S_{\rm bd}(\mbf h^{(q)},\signvec;N)$. Note that from \eqref{definition for inner points}, the boundary points that violate this definition have at least one coordinate inside the intervals of the form (for some $A$)
$$[A, (1+\delta^2_N)A] \text{ or } [(1+\delta_N)(1-\delta_N^2)A, (1+\delta_N)A],$$ 
which are of length $O(\delta^2_N A) $, while the original interval is of length $\delta_NA$. Thus the boundary is a finite union of $O_d(1)$ boxes whose expected cardinality is
\[
\ll_d\frac{\delta^2_N}{\delta_N}E(\mbf h^{(q)},\signvec;N)\ll\frac{E(\mbf h^{(q)},\signvec;N)}{(\log N)^2}.
\]
As in Lemmas \ref{Proj two number of elements in S for N^d/4} and \ref{Proj two d_N error estimation for the number of elements}, each such boundary box is a finite signed sum of anchored boxes, and the number of anchored boxes needed is bounded only in terms of $d$. Hence Beck’s discrepancy estimate \eqref{Thm of Beck} applies uniformly.” Therefore, the discrepancy error is $O((\log N)^{d+1/2})$, which is also $O(E(\mbf h^{(q)},\signvec;N)/(\log N)^2)$ because
\[
E(\mbf h^{(q)},\signvec;N)\gg \delta_N^{d+1}(\log N)^{3d+6}=(\log N)^{d+4}.
\]
Thus we have
\beq\label{Proj two boundary estimate revised}
|S_{\rm bd}(\mbf h^{(q)},\signvec;N)|\ll_{\bma}\frac{|S(\mbf h^{(q)},\signvec;N)|}{(\log N)^2}.
\eeq

Combining \eqref{Proj two AP estimate revised} and \eqref{Proj two boundary estimate revised}, first with $\signvec=\signvec^+$ and then with $\signvec=\signvec^-$, gives
\[
\begin{aligned}
&\left|\sum_{S(\mbf h^{(q)},\signvec^+;N)}e^{2\pi \i\L(\mbf n)}-\sum_{S(\mbf h^{(q)},\signvec^-;N)}e^{2\pi \i\L(\mbf n)}\right|\\
&\qquad\ll_{\bma}\frac{|S(\mbf h^{(q)},\signvec^+;N)|+|S(\mbf h^{(q)},\signvec^-;N)|}{(\log N)^2},
\end{aligned}
\]
which contradicts the assumption that $\mbf h^{(q)}$ is $\mathfrak p$-big. Hence each special line contains at most one $\mathfrak p$-big vector.
\end{proof}

\begin{proof}[Proof of Lemma \ref{Proj two number of e-big vectors}]
It remains to count special lines. Use the coordinates
\[
u_0=l_1,\qquad
u_j=l_1-l_{j+1}\quad(1\le j\le d-1),\qquad
u_d=l_{d+1}.
\]
By \eqref{Proj two condition for l}, each $u_i$ lies in an interval of length $O(\log N/\delta_N)$ and has lower admissible boundary $\gg\log\log N/\delta_N$. Under one step along a special line,
\[
u_0\mapsto u_0+A_N,\qquad
u_j\mapsto u_j+2A_N\quad(1\le j\le d-1),\qquad
u_d\mapsto u_d+(d+1)A_N.
\]
A maximal special line is determined by its first point. For this first point, at least one of the coordinates $u_i$ must lie within $O(A_N)$ of its lower admissible boundary; otherwise the point could be shifted one step backwards and would remain in $L_2(N)$. Therefore the number of possible first points is
\[
\ll A_N\left(\frac{\log N}{\delta_N}\right)^d
\ll(\log\log N)(\log N)^d\delta_N^{-(d+1)},
\]
because $A_N\ll\log\log N/\delta_N$. By Lemma \ref{Proj two key lemma}, each special line contains at most one $\mathfrak p$-big vector. This proves Lemma \ref{Proj two number of e-big vectors}.
\end{proof}

\normalsize

\bibliographystyle{plainnat} 
\bibliography{references_revision3}

\end{document}